\documentclass[11pt,leqno]{article}
\usepackage{amsmath, amsfonts, amssymb, amsthm, graphicx}

\numberwithin{equation}{section}

\newtheorem{theorem}{Theorem}
\newtheorem{lemma}{Lemma}
\newtheorem{remark}{Remark}

\newtheorem{corollary}{Corollary}

\usepackage[dvips]{epsfig}
\usepackage{graphicx}

\begin{document}

\title{\textbf{Wave breaking in the Ostrovsky--Hunter equation}}
\author{Yue Liu$^1$, Dmitry Pelinovsky$^2$,  and Anton Sakovich$^2$ \\
{\small $^{1}$ Department of Mathematics, University of Texas at Arlington, Arlington, TX, 76019, USA} \\
{\small $^{2}$ Department of Mathematics, McMaster
University, Hamilton, ON, L8S 4K1, Canada}}
\date{}
\maketitle

\begin{abstract}
The Ostrovsky--Hunter equation governs evolution of shallow water
waves on a rotating fluid in the limit of small high-frequency
dispersion. Sufficient conditions for the wave breaking in the
Ostrovsky--Hunter equation are found both on an infinite line and
in a periodic domain. Using the method of characteristics, we also
specify the blow-up rate at which the waves break. Numerical
illustrations of the finite-time wave breaking are given in a periodic domain.
\end{abstract}

\section{Introduction}

The nonlinear evolution equation
\begin{equation}
\label{Ostrovsky} \left( u_t + u u_x - \beta u_{xxx} \right)_x = \gamma u,
\end{equation}
with $\beta, \gamma \in \mathbb{R}$ and $u(t,x) : \mathbb{R}_+
\times \mathbb{R} \mapsto \mathbb{R}$, was derived by Ostrovsky
\cite{Ostrov} to model small-amplitude long waves in a rotating
fluid of a finite depth. This equation generalizes the
Korteweg--de Vries equation (that corresponds to $\gamma = 0$) by
the additional term induced by the Coriolis force. Mathematical
properties of the Ostrovsky equation (\ref{Ostrovsky}) were
studied recently in many details, including the local and global
well-posedness in energy space \cite{GL07,LM06,Ts,VL04}, stability of
solitary waves \cite{LL06,L07,VL05}, and convergence of solutions in
the limit of the Korteweg--de Vries equation \cite{LL07,VL05}.

We shall consider the limit of no high-frequency
dispersion $\beta = 0$, when the evolution equation
(\ref{Ostrovsky}) can be written in the form
\begin{equation}
\label{OstHunter}
\left\{ \begin{array}{l} (u_t + u u_x)_x = \gamma u, \quad t > 0,\\
u(0,x) = u_0(x), \end{array} \right.
\end{equation}
where $x$ is considered either on a circle or on an infinite line.
In this form, the main equation (\ref{OstHunter}) is known under
different names such as the reduced Ostrovsky equation \cite{Parkes,Step}, the
Ostrovsky--Hunter equation \cite{Bo}, the short-wave equation
\cite{Hu}, and the Vakhnenko equation \cite{Vakh1,Vakh2}. We shall use
the name of the Ostrovsky--Hunter equation for convenience. We
also consider $\gamma > 0$ since the other case $\gamma < 0$ is
covered by the reflection $x \to -x$ and $u \to -u$ of the
solutions for $\gamma > 0$.

According to the method of characteristics, the simple-wave
equation $u_t + u u_x = 0$ (that corresponds to $\gamma = 0$)
develops wave breaking in a finite time for any initial data
$u(0,x) = u_0(x)$ on an infinite line or in a periodic domain if
$u_0(x)$ is continuously differentiable and there is a point $x_0$
such that $u_0'(x_0) < 0$. More precisely, we say that the
finite-time wave breaking occurs if there exists a finite time $T
\in (0,\infty)$ such that
\begin{equation}
\label{wave-breaking-criterion}
\liminf_{t \uparrow T} \inf_x u_x(t,x) = -\infty, \quad \mbox{while}
\quad \limsup_{t \uparrow T} \sup_x \left| u(t,x) \right| < \infty.
\end{equation}
In view of the result for $\gamma = 0$, we address the question if the
low-frequency dispersion term with $\gamma > 0$ in the
Ostrovsky--Hunter equation (\ref{OstHunter}) can stabilize global
dynamics of the simple-wave equation $u_t + uu_x = 0$.

Hunter \cite{Hu} found a sufficient condition for wave breaking of
the Cauchy problem (\ref{OstHunter}) in a periodic domain and
provided numerical evidences of the finite-time wave breaking for
the sinusoidal initial data $u_0(x)$. To be precise, the
main result of \cite{Hu} can be formulated as follows.

\begin{theorem}[Hunter, 1990]
Let $u_0(x) \in C^1(\mathbb{S})$, where $\mathbb{S}$ is a circle
of unit length and define
$$
\inf_{x \in \mathbb{S}} u_0'(x) = -m \quad \mbox{\rm and}
\quad \sup_{x \in \mathbb{S}} \left| u_0(x) \right| = M.
$$
If $m^3 > 4 M (4 + m)$, a smooth solution $u(t,x)$ of the Cauchy problem
(\ref{OstHunter}) with $\gamma = 1$ breaks down at a finite time $T \in (0,2
m^{-1})$. \label{theorem-hunter}
\end{theorem}

We shall study wave breaking of the Cauchy problem (\ref{OstHunter}) in more details.
First, the Cauchy problem is shown to be locally well-posed if $u_0 \in H^s$, $s \geq 2$
both on an infinite line $\mathbb{R}$ and on a unit circle $\mathbb{S}$. The blow-up
alternative is derived to claim that the solution $u(t,x)$ blows up in a finite
time in the $H^s$ norm with $s \geq 2$ if and only if $\inf_x u_x(t,x)$ becomes unbounded
from below. Using the integral estimates and the method of characteristics similar to
analysis of the Camassa--Holmes equation in \cite{Constantin,C-E} and
the Degasperis--Processi equation in \cite{LY1,LY2}, we find various
sufficient conditions for the wave breaking, which are sharper than Theorem \ref{theorem-hunter}.
Moreover, we also obtain the blow-up rate at which the waves break in a finite time.

We note that, unlike the Ostrovsky equation (\ref{Ostrovsky}), the
Ostrovsky--Hunter equation (\ref{OstHunter}) is integrable using
the inverse scattering transform method \cite{Vakh2}. This method
allows us to solve the initial-value problem (\ref{OstHunter})
formally by working with the spectral theory for a third-order
differential operator, which is similar to the Lax operator for
the Hirota--Satsuma equation \cite{Satsuma}. As a particular
property of an integrable model, the Ostrovsky--Hunter equation
(\ref{OstHunter}) has a hierarchy of conserved quantities, which
follows from results of \cite{Vakh3} and \cite{Satsuma} after
exchanging densities and fluxes. This hierarchy includes the first
two conserved quantities
\begin{equation}
Q = \int u^2 dx, \quad E = \int \left[ \gamma (\partial_x^{-1} u)^2 + \frac{1}{3} u^3 \right] dx,
\end{equation}
where the anti-derivative operator is defined by the integration of $u(t,x)$ in $x$ subject to the zero-mass
constraint $\int u dx = 0$. Higher-order conserved quantities of the Ostrovsky--Hunter equation (\ref{OstHunter})
involve higher-order anti-derivatives, which are defined under additional constraints on the solution $u(t,x)$.
Therefore, conserved quantities of the Ostrovsky--Hunter equation are not related to the $H^s$-norms
and hence are not so useful in the study of global well-posedness in the energy space,
in a sharp contrast with very similar short-pulse and Hirota--Satsuma equations
studied in \cite{PelSak} and \cite{iorio}, respectively. Our analysis does not rely, therefore, on
integrability properties of the Ostrovsky--Hunter equation (\ref{OstHunter}). We also
emphasize that integrability of the nonlinear evolution equations does not prevent a finite-time blow-up,
see \cite{Constantin,C-E,LY1,LY2} for analysis of wave breaking in other integrable equations.

The other problem related to the subject of this paper is the existence and stability of
spatially periodic and localized traveling-wave solutions $u(t,x) = \varphi(x-ct)$ with speed $c \in \mathbb{R}$.
Function $\varphi(x)$ is defined by solutions of the differential equation
\begin{equation}
\label{stationary-wave}
(c-\varphi(x)) \varphi''(x) = (\varphi'(x))^2 - \gamma \varphi(x),
\end{equation}
where $x$ is considered either on a circle or on an infinite line.
Bounded $2\pi$-periodic solutions $\varphi(x)$ were shown in
\cite{Bo} to exist for the wave speeds
$$
1 \leq \frac{c}{\gamma} \leq \frac{\pi^2}{9},
$$
where $c = \gamma$ corresponds to the small-amplitude sinusoidal
wave and $c = \frac{\pi^2}{9}\gamma $ corresponds to the
large-amplitude crest wave (called the parabolic wave) which is
given by the piecewise continuous quadratic polynomial in $x$
$$
\varphi(x) = \frac{\gamma}{16} (3 x^2 - \pi^2), \quad x \in
[-\pi,\pi],
$$
with a discontinuous slope at the crests located at $x = \pm \pi$.
Analytical and numerical approximations of the periodic wave
solutions can be found in \cite{Bo} and \cite{Hu}. Our results on
wave breaking in a periodic domain $\mathbb{S}$ do not exclude
possibility of global well-posedness of the Cauchy problem for
small initial data $u_0(x)$ and stability of periodic wave solutions
satisfying (\ref{stationary-wave}). The latter problems are left, however, beyond the scopes of this
paper.

No localized solutions $\varphi(x)$ were found on a real line $\mathbb{R}$, except for multi-valued
loop solitons \cite{Vakh1} and other exotic solutions \cite{Parkes,Step} that do not belong to
$H^s(\mathbb{R})$ with $s \geq 2$. In Appendix A, we prove that no classical
solutions $\varphi(x) \in C^2(\mathbb{R})$ with decay
$$
\lim_{|x| \to \infty} \varphi(x) = \lim_{|x| \to \infty} \varphi'(x) = 0
$$
exist. Again, global well-posedness of the Cauchy problem on an infinite line
for small initial data $u_0(x)$ is not excluded by our results.

The paper is organized as follows. Section 2 gives a sufficient condition for the wave
breaking in a periodic domain. The blow-up rate of the wave breaking is studied in Section 3
based on the method of characteristics. Similar results on an infinite line are reported in Section 4.
Section 5 gives numerical evidences of wave breaking in a periodic domain. Appendix A contains
results on non-existence of localized traveling-wave solutions.

\section{Wave breaking in a periodic domain}

Let $\gamma > 0$ and $\mathbb{S}$ denote a circle of a unit length.
Local well-posedness of the Cauchy problem (\ref{OstHunter}) with
initial data $u_{0}\in H^{s}(\mathbb{S})$, $s \ge 2$ can be
obtained using the work of Sch\"{a}fer \& Wayne \cite{ScWa}
who studied a very similar short-pulse equation
\begin{equation}
\label{short-pulse} \left\{ \begin{array}{l} u_{xt} = u + (u^3)_{xx}, \quad t > 0,\\
u(0,x) = u_0(x), \end{array} \right.
\end{equation}
on an infinite line. More precisely, we have the following local
well-posedness result.

\begin{lemma}  Assume that $u_{0}(x) \in H^{s}(\mathbb{S})$, $s \ge 2$
and $\int_{\mathbb{S}} u_0(x) \, dx = 0$. Then  there exist a
maximal time $T=T(u_{0})>0$  and a unique solution $u(t,x)$ to the
Cauchy problem (\ref{OstHunter}) such that
$$
u(t,x) \in C([0,T);H^{s}(\mathbb{S}))\cap
C^{1}([0,T);H^{s-1}(\mathbb{S}))
$$
with the following two conserved quantities
$$
\int_{\mathbb{S}} u(t, x) dx  = 0,  \qquad t \in [0, T)
$$
and
$$
\int_{\mathbb{S}} u^2(t, x) dx = \int_{\mathbb{S}} u_0^2(x) dx,    \qquad  t \in [0, T).
$$
Moreover, the solution depends continuously on the initial data,
i.e. the mapping $u_{0}\mapsto u : \; H^{s}(\mathbb{S})
\rightarrow C([0,T);H^{s}(\mathbb{S}))\cap
C^{1}([0,T);H^{s-1}(\mathbb{S}))$ is continuous.
\label{lemma-wellposedness}
\end{lemma}

\begin{proof}
Existence, uniqueness, and continuous dependence in
$H^s(\mathbb{R})$, $s \geq 2$ were proven in \cite{ScWa} in the
content of the short-pulse equation (\ref{short-pulse}). The same
method based on modified Picard's iterations works in
$H^s(\mathbb{S})$, so that the first part of Lemma
\ref{lemma-wellposedness} is an extension of the main theorem of
\cite{ScWa} to a periodic domain. To prove the zero-mass
constraint, we note that $u_t \in C([0, T);H^1)$ and $uu_x \in
C([0, T);H^1)$ for the solution $u \in C([0,T);H^{s}(\mathbb{S}))
\cap C^{1}([0,T);H^{s-1}(\mathbb{S}))$ with $s \geq 2$. By
Sobolev's embedding of $H^1(\mathbb{S})$ to $C(\mathbb{S})$, we
obtain
$$
\gamma \int_{\mathbb{S}} u(t, x) \, dx =
\int_{\mathbb{S}} u_{tx} \, dx + \int_{\mathbb{S}} ( u u_x)_x \, dx = 0, \quad
t \in (0,T).
$$
To prove conservation of the $L^2$-norm, we consider the balance
equation
$$
(u^2)_t = \left( \gamma (\partial_x^{-1} u)^2 - \frac{2}{3} u^3
\right)_x, \quad x \in \mathbb{S}, \quad t \in (0,T),
$$
where $\partial_x^{-1} u = \gamma^{-1} (u_t + u u_x) \in C([0,T);H^1(\mathbb{S}))$, so that
$\partial_x^{-1} u(t,x)$ is continuous on $\mathbb{S}$ for all $t \in [0,T)$.
Integrating the balance equation, we obtain
$$
\int_{\mathbb{S}} u^2(t, x) \, dx = \int_{\mathbb{S}} u_0^2(x) \, dx, \quad
t \in [0,T).
$$
This completes the proof of Lemma \ref{lemma-wellposedness}.
\end{proof}

\begin{remark}
The assumption $\int_{\mathbb{S}} u_0(x) \, dx = 0 $ in Lemma \ref{lemma-wellposedness}
on the initial data $u_0$ is necessary. In fact, the zero-mass constraint on $u_0$
can be derived from the following estimate
$$
\left | \int_{\mathbb{S}} u(t, x) \, dx - \int_{\mathbb{S}} u_0 (x) \, dx \right | \le \| u(t, \cdot ) - u_0(\cdot) \|_{L^2(\mathbb{S})}, \quad \forall t \in (0, T).
$$
Note that $ \int_{\mathbb{S}} u(t, x) \, dx = 0, $ for all $ t \in
(0, T) $ and $ u \in C([0, T); H^2(\mathbb{S})).$ Hence the above
estimate implies in the limit $t \downarrow 0$ that $ \int_{\mathbb{S}} u_0(x) dx = 0$. Note
that no zero-mass constraint on $u_0$ is necessary on an infinite line
\cite{ScWa}.
\end{remark}

\begin{remark}
\label{remark-independence}
The maximal time $T > 0$ in Lemma \ref{lemma-wellposedness} is independent of $s \geq 2$
in the following sense. If $u_0(x) \in H^s(\mathbb{S}) \cap H^{s'}(\mathbb{S})$
for $s,s' \geq 2$ and $s \neq s'$, then
\begin{eqnarray*}
u(t,x) \in C([0,T);H^{s}(\mathbb{S}))\cap C^1([0,T);H^{s-1}(\mathbb{S}))
\end{eqnarray*}
and
\begin{eqnarray*}
u(t,x) \in C([0,T');H^{s'}(\mathbb{S}))\cap C^1([0,T');H^{s'-1}(\mathbb{S}))
\end{eqnarray*}
with the same $T' = T$. See Yin \cite{Yin} for an adaptation of the Kato method
\cite{Kato} to the proof of this statement.
\end{remark}

By using the local well-posedness result in Lemma
\ref{lemma-wellposedness} and energy estimates, one can get the
following precise blow-up scenario of the solutions to the Cauchy
problem (\ref{OstHunter}).

\begin{lemma}
Let $u_0(x) \in H^s(\mathbb{S})$, $s \geq 2$ and $u(t,x)$ be a
solution of the Cauchy problem (\ref{OstHunter}) in Lemma
\ref{lemma-wellposedness}. The solution blows up in a finite time
$T \in (0,\infty)$ in the sense of $\lim_{t \uparrow T} \|
u(t,\cdot)\|_{H^s} = \infty$ if and only if
$$
\lim_{t \uparrow T} \inf_{x \in \mathbb{S}} u_{x}(t,x) = -
\infty.
$$
\label{lemma-blowup}
\end{lemma}

\begin{proof}
Assume a finite maximal existence time $T \in (0,\infty)$ and
suppose there is $M > 0$ such that
\begin{equation}
\label{bound-der} \inf_{x \in \mathbb{S}} u_{x}(t,x)  \ge -M, \quad \forall t \in [0,T).
\end{equation}
Applying density arguments, we approximate initial value
$u_0(x) \in H^2(\mathbb{S})$ by functions $u_0^n(x) \in H^3(\mathbb{S})$, $n
\geq 1$, so that $\lim_{n \to \infty} u_0^n = u_0$. Furthermore,
write $u^n(t,x)$ for the solution of the Cauchy problem
(\ref{OstHunter}) with initial data $u_0^n(x) \in H^3(\mathbb{S})$.
Using the regularity
result proved in Lemma \ref{lemma-wellposedness}, it follows from
Sobolev's embedding that, if $u^n(t,x) \in C([0,T);H^3(\mathbb{S}))$,
then $u^n(t,x)$ is a twice continuously differentiable periodic
function of $x$ on $\mathbb{S}$ for any $t \in [0,T)$.  It is then
deduced from the Ostrovsky--Hunter equation (\ref{OstHunter}) that
$$
\frac{d}{dt} \int_{\mathbb{S}}(u^n_x)^2 dx = - \int_{\mathbb{S}} u^n_x
(u^n_x)^2 dx \le M \int_{\mathbb{S}} (u^n_x)^2 dx,
$$ and
$$
\frac{d}{dt} \int_{\mathbb{S}} (u^n_{xx})^2 dx = - 5 \int_{\mathbb{S}}
u^n_x (u^n_{xx})^2 dx \le 5 M
\int_{\mathbb{S}} (u^n_{xx})^2 dx,
$$
where we have used the uniform bound (\ref{bound-der}). The
Gronwall inequality then yields
$$
\|u^n_x \|_{L^2} \le \| (u_0^n)' \|_{L^2} e^{\frac{M}{2} t}, \quad
$$ and
$$
\|u^n_{xx}\|_{L^2} \le \| (u_0^n)'' \|_{L^2} e^{\frac{5}{2} Mt},
\quad 0 \le t < T.
$$
Since $ \|u^n_0\|_{H^2} $ converges to $\|u_0\|_{H^2}$ as $n \to
\infty$, we infer from  the continuous dependence
of the local solution $u$ on initial data $u_0$ that the norm in
$H^2(\mathbb{S})$ of the solution $u$ in Lemma
\ref{lemma-wellposedness} does not blow up in the finite time $T <
\infty$ and therefore either $T$ is not a maximal existence time
or the bound (\ref{bound-der}) does not hold as $t \uparrow T$.
Since $T$ is independent on $s \geq 2$ by Remark \ref{remark-independence}, the norm in
$H^s(\mathbb{S})$ for any $s \geq 2$ of the solution in Lemma
\ref{lemma-wellposedness} blows up in a finite time $T \in
(0,\infty)$ if and only if the bound (\ref{bound-der}) does not
hold as $t \uparrow T$.
\end{proof}

The main result of this section is the following sufficient
condition for the wave breaking in the Ostrovsky--Hunter equation
(\ref{OstHunter}).

\begin{theorem}
\label{theorem-wave-breaking} Assume that $u_{0}(x) \in
H^{s}(\mathbb{S})$, $s \ge 2$ and $\int_{\mathbb{S}} u_0(x) \, dx
= 0$. If $u_0$ satisfies either
\begin{equation}
\label{condition-one}
\int_{\mathbb{S}} \left( u_0'(x) \right)^3 \, dx < - \left(
\frac{3 \gamma}{2} \| u_0 \|_{L^2} \right)^{3/2}
\end{equation}
or
\begin{equation}
\label{condition-two}
\int_{\mathbb{S}} \left( u_0'(x) \right)^3 \, dx < 0 \quad
\mbox{and} \quad \| u_0 \|_{L^2} > \frac {3 \gamma}{4},
\end{equation}
then the solution $u(t,x)$ of the Cauchy problem (\ref{OstHunter})
in Lemma \ref{lemma-wellposedness}
blows up in finite time $T \in (0,\infty)$ in the sense of Lemma \ref{lemma-blowup}.
\end{theorem}

\begin{proof}
Let $T > 0$ be the maximal time of existence of the solution $u(t,x)$
in Lemma \ref{lemma-wellposedness}.  Then, we obtain the a priori
differential estimate
\begin{eqnarray*}
\frac{d}{dt}\int_{\mathbb{S}} u_{x}^3\,dx & = & 3
\int_{\mathbb{S}} u_x^2 \left( - u_x^2 - u u_{xx} + \gamma u
\right) \,dx \\ & = & -2 \int_{\mathbb{S}} u_x^4 \,dx
+ 3 \gamma \int_{\mathbb{S}} u u_x^2 \,dx \\
& \leq &  -2 \| u_x \|^4_{L^4} + 3 \gamma \| u \|_{L^2} \| u_x \|^2_{L^4}\\
& = &  -2 \left( \| u_x \|^2_{L^4} - \frac{3 \gamma}{4} \| u_0 \|_{L^2}
\right)^2 + \frac{9 \gamma^2}{8} \| u_0 \|_{L^2}^2,
\end{eqnarray*}
where we have used the Cauchy--Schwarz inequality and the $L^2$-norm conservation.
An application of H\"{o}lder's inequality yields
\begin{equation}
\label{Holder-inequality}
\| u_x \|^3_{L^3} \leq \| u_x \|^3_{L^4}.
\end{equation}
Let $V(t) = \int_{\mathbb{S}}u_x^3(t,x)\,dx$ for all $t \in [0,T)$,
$Q_0 = \| u_0 \|_{L^2}$, and assume that
$$
V(0) < -\left( \frac{3 \gamma Q_0}{2} \right)^{\frac{3}{2}} < 0.
$$
Then, we have
$$
\| u_x \|^2_{L^4} - \frac{3 \gamma}{4} \| u_0 \|_{L^2} \geq
\| u_x \|^2_{L^3} - \frac{3 \gamma}{4} \| u_0 \|_{L^2} \geq |V|^{\frac{2}{3}} - \frac{3 \gamma Q_0}{4},
$$
so that the a priori differential inequality is closed at
$$
\frac{dV}{dt} \leq -2 \left( |V|^{\frac{2}{3}} -\frac{3 \gamma Q_0}{4}
\right)^2 + \frac{9 \gamma^2 Q_0^2}{8},
$$
where the right-hand-side is negative at $t = 0$. By the continuation
argument, $V(t)$ is decreasing on $[0,T)$ so that $V(t) \leq V(0) < 0$.
We need to prove that $T$ is finite and
$\lim_{t \uparrow T} V(t) = -\infty$. Let $y =
|V|^{1/3}$ and obtain that
$$
\frac{dy}{dt} \geq \frac{2}{3} \left( y^2 - \frac{3 \gamma Q_0}{2}
\right),
$$
where the right-hand-side is positive at $t = 0$. By the comparison
principle for differential equations $y(t) \geq y^+(t)$ for all $t
\in [0,T)$, where $y^+(t)$ solves the differential equation
$$
\left\{ \begin{array}{l} \dot{y}^+ = \frac{2}{3} \left( (y^+)^2 -
\frac{3 \gamma Q_0}{2} \right), \\
y^+(0) = y(0)
\end{array} \right.
$$
Since $y(0) > \left( \frac{3 \gamma Q_0}{2} \right)^{\frac{1}{2}}$,
there is a finite time $T^+ \in (0,\infty)$ such that $\lim_{t
\uparrow T^+} y^+(t) = +\infty$ and therefore, there is a time $T
\in (0,T^+)$ such that $\lim_{t \uparrow T^+} y(t) = +\infty$.

To prove the second sufficient condition (\ref{condition-two}), we note that since
$\int_{\mathbb{S}} u(t,x) dx = 0$,  for each $t\in [0,T)$ there is a
$\xi_{t}\in \mathbb{S}$ such that $u(t,\xi_{t})=0$. Then for any $x \in  [\xi_t,\xi_t+\frac{1}{2}]$,
we have
\begin{equation*}
u^{2}(t,x)=\left( \int_{\xi_t}^{x}u_x(t,x) \,dx\right)^2\leq
(x-\xi_t)\int^x_{\xi_t}u_x^2(t,x)\,dx \leq \frac{1}{2} \| u_x \|_{L^2}^2.
\end{equation*}
Combining it with a similar estimate on $[\xi_t+\frac{1}{2},\xi_t +1]$
thanks to periodicity of $u(t,x)$ in $x$ for all $t \in [0,T)$, we have
\begin{equation*}
\sup_{x\in \mathbb{S}}u^{2}(t,x)\leq
\frac{1}{2} \| u_x \|^2_{L^2} \leq \frac {1}{2} \| u_x \|^2_{L^4}.
\end{equation*}
Therefore, continuing the a priori differential inequality above, we
obtain
\begin{eqnarray*}
\frac{d}{dt}\int_{\mathbb{S}} u_{x}^3\,dx & \leq &
-2 \| u_x \|^4_{L^4} + 3 \gamma \| u \|_{L^2} \| u_x \|^2_{L^4} \\
& \leq &
-2 \| u_x \|^4_{L^4} + 3 \gamma  \| u_0 \|_{L^2}^{-1}
\| u \|_{L^{\infty}}^2 \| u_x \|^2_{L^4} \\
& \leq & - \alpha  \| u_x \|_{L^4}^4,
\end{eqnarray*}
where $\alpha =  2 - \frac {3 \gamma }{2 \| u_0 \|_{L^2}} > 0 $ by
the assumption. By the same H\"{o}lder's inequality, we obtain
\begin{equation*}
\frac{d V}{dt} \leq -\alpha |V|^{\frac{4}{3}}.
\end{equation*}
where $V(t)=\int_{\mathbb{S}}u_x^3(t,x)\,dx$ for all $t \in
[0,T)$ and $V(0) < 0$ is assumed. Then, by the comparison
principle, $V(t) \leq V^-(t)$ for all $t \in [0,T)$, where
$V^-(t)$ solves the differential equation
$$
\left\{ \begin{array}{l} \dot{V}^- = - \alpha (V^-)^{\frac{4}{3}},
\\ V^-(0) = V(0). \end{array} \right.
$$
Since $V(0) < 0$, there is a finite time $T^- \in (0,\infty)$ such that
$\lim_{t \uparrow T} V^-(t) = -\infty$ and therefore,
there is a finite time $T \in (0,T^-)$ such that $\lim_{t \uparrow
T} V(t) = - \infty$. In both cases, we have
$$
\inf_{x \in \mathbb{S}} u_x^3(t,x) \le \int_{\mathbb{S}} u_x^3 \,
dx \equiv V(t),
$$
which implies immediately that
$$
\lim_{t \uparrow T} \inf_{x \in \mathbb{S}} u_{x} (t,x) = -
\infty.
$$
This completes the proof of the theorem.
\end{proof}

\begin{remark}
Let $\inf_{x \in \mathbb{S}} u_0'(x) = -m$. The first sufficient condition
(\ref{condition-one}) in Theorem \ref{theorem-wave-breaking}
can be rewritten as
$$
m^2 > \frac{3 \gamma}{2} \| u_0 \|_{L^2},
$$
which reminds us the sufficient condition in Theorem \ref{theorem-hunter}
for $\gamma = 1$ given by
$$
m(m^2 - 4 \| u_0 \|_{L^{\infty}}) > 16 \| u_0 \|_{L^{\infty}}.
$$
If $\|u_0 \|_{L^2}$ (and then $\| u_0 \|_{L^{\infty}}$) is large, the
slope of $u_0'(x)$ has to be steep enough to lead to the wave breaking.
In a contrast, the second sufficient condition (\ref{condition-two})
in Theorem \ref{theorem-wave-breaking} shows that  any smooth initial
profile with $\int_{\mathbb{S}} (u_0'(x))^3 dx < 0$
and sufficiently large $\| u_0 \|_{L^2}$ breaks in a finite time.
\end{remark}

\section{Blow-up rate of wave breaking}

We shall investigate here the blow-up rate of the wave breaking
for solutions of the Cauchy problem (\ref{OstHunter}), which we rewrite here as
\begin{equation}
\label{cauchy}
\left\{ \begin{array}{l} u_t + u u_x = \gamma \partial_x^{-1} u, \quad  x \in \mathbb{S}, \;\;
t > 0, \\
u(0,x) = u_0(x), \qquad \phantom{t} x \in \mathbb{S}, \end{array} \right.
\end{equation}
where $\partial_x^{-1}$ is the mean-zero anti-derivative in the sense of
\begin{equation}
\label{anti-derivative}
\partial_x^{-1} u = \int_0^x u(t,x') dx' - \int_{\mathbb{S}} \int_0^x u(t,x') dx' dx.
\end{equation}
We use the method of characteristics, which is also used in a similar context
by Hunter \cite{Hu}.
Let $T > 0$ be the maximal time of existence of the solution $u(t,x)$ of the Cauchy
problem (\ref{cauchy}) in Lemma \ref{lemma-wellposedness} with the initial data
$u_0 \in H^s(\mathbb{S})$ for $s \geq 2$. For all $t \in [0,T)$ and $\xi \in \mathbb{S}$,
define
$$
x = X(t,\xi), \quad u(t,x) = U(t,\xi), \quad \partial_x^{-1} u(t,x) = G(t,\xi),
$$
so that
\begin{equation}
\label{characteristics}
\left\{ \begin{array}{l} \dot{X}(t) = U, \\ X(0) = \xi, \end{array} \right. \quad
\left\{ \begin{array}{l} \dot{U}(t) = \gamma G, \\ U(0) = u_0(\xi), \end{array} \right.
\end{equation}
where dots denote derivatives with respect to time $t$ on a particular characteristics $x = X(t,\xi)$
for a fixed $\xi \in \mathbb{S}$. Applying classical results in the theory of ordinary differential
equations, we obtain the following two useful results on the solutions of the initial-value problem
(\ref{characteristics}).

\begin{lemma}
Let $u_{0}(x) \in H^{s}(\mathbb{S})$, $s \geq 2$ and $T>0$ be the
maximal existence time of the solution $u(t,x)$ in Lemma \ref{lemma-wellposedness}.
Then there exists a unique solution $X(t,\xi) \in C^{1}([0,T)\times \mathbb{S})$ to
the initial-value problem (\ref{characteristics}). Moreover, the map
$X(t,\cdot) : \mathbb{S} \mapsto \mathbb{R}$ is an increasing diffeomorphism with
$$
\partial_{\xi} X(t,\xi) = \exp \left ( \int_{0}^{t} u_x(s,X(s,\xi)) ds \right ) > 0,\, \;
\forall t \in [0,T), \;\; \forall x \in \mathbb{S}.
$$
\label{lemma-characteristics}
\end{lemma}

\begin{proof}
Consider the integral equation
$$
X(t,\xi) = \xi + \int_0^t u(s,X(s,\xi)) ds, \quad t \in [0,T),
$$
where $u(t,x) \in C([0,T);H^{s}(\mathbb{S}))\cap
C^{1}([0,T);H^{s-1}(\mathbb{S}))$ for $s \geq 2$, according to Lemma \ref{lemma-wellposedness}.
By the ODE theory, there exists a unique solution $X(t,\xi) \in C^{1}([0,T)\times \mathbb{S})$
of the integral equation above. Using the chain rule, we obtain
$$
\partial_{\xi} \dot{X} = u_x(t,X(t,\xi)) \partial_{\xi} X  \quad \Rightarrow \quad
\partial_{\xi} X(t,\xi) = \exp\left(\int_0^t u_x(t,X(s,\xi)) ds\right),
$$
so that $\partial_{\xi} X(t,\xi) > 0$ for all $t \in [0,T)$ and $\xi \in \mathbb{S}$.
\end{proof}

\begin{lemma}
\label{lemma-bounds}
Let $u_{0}(x) \in H^{s}(\mathbb{S})$, $s \ge 2$ and $T>0$ be the
maximal existence time of the solution $u(t,x)$ in Lemma \ref{lemma-wellposedness}. Then
the solution $u(t,x)$ satisfies
$$
\sup_{s \in [0,t]} \| u(s, \cdot)\|_{L^{\infty}} \leq \| u_{0}\|_{L^{\infty}} + \gamma t \| u_{0}\|_{L^{2}},
\quad \forall t\in [0,T).
$$
\end{lemma}

\begin{proof}
By Lemma \ref{lemma-characteristics}, the function $x = X(t,\xi)$ is invertible
in $\xi \in \mathbb{S}$ for any $t \in [0,T)$. Then, we have
$$
\sup_{s \in [0,t]} \sup_{x \in \mathbb{S}} |u(s,x) | = \sup_{s \in [0,t]} \sup_{\xi \in \mathbb{S}} |U(s,\xi)|,
\quad t \in [0,T).
$$
Since $\partial_x^{-1} u(t,x) \in C([0,T);H^{s+1}(\mathbb{S}))$ is the mean-zero periodic function of $x$
for each $t \in [0,T)$, there exists a $\xi_t \in \mathbb{S}$ such that
$\partial_x^{-1} u(t,\xi_t) = 0$. Then for any $x\in \mathbb{S}$ and $t \in [0,T)$, we have
\begin{equation*}
|\partial_x^{-1} u(t,x)| \leq \left| \int_{\xi_t}^{x} u(t,x) \,dx \right| \leq
\int_{\mathbb{S}} |u(t,x)| dx \leq \| u_0 \|_{L^2},
\end{equation*}
where we use the Cauchy--Schwarz inequality and the $L^2$ norm conservation.
Using the integral equation
$$
U(t,\xi) = u_0(\xi) + \gamma \int_0^t G(s,\xi) ds, \quad t \in [0,T),
$$
we obtain
\begin{eqnarray*}
\sup_{s \in [0,t]} \sup_{x \in \mathbb{S}} |u(s,x) | & \leq & \| u_0 \|_{L^{\infty}} + \gamma t \sup_{s \in [0,t]} \sup_{\xi \in \mathbb{S}} |G(s,\xi)| \\
& = & \| u_0 \|_{L^{\infty}} + \gamma t \sup_{s \in [0,t]} \sup_{x \in \mathbb{S}}
\left| \partial_x^{-1} u(s,x) \right| \\
& \leq & \| u_{0}\|_{L^{\infty}} + \gamma t  \| u_{0}\|_{L^{2}}, \quad t \in [0,T),
\end{eqnarray*}
and the lemma is proved.
\end{proof}

Using the method of characteristics, we obtain a sufficient condition for
the wave breaking in the Cauchy problem (\ref{cauchy}) that is different from
the sufficient conditions of Theorem \ref{theorem-wave-breaking}.

\begin{theorem}
\label{theorem-blowup-characteristics} Let $\varepsilon>0$ and
$u_{0}(x) \in H^{s}(\mathbb{S})$, $s \ge 2$. Let $T_1$ be the
smallest positive root of
\begin{equation}
\label{time-bound}
2 \sqrt{\gamma} T_1 \left( \| u_0 \|_{L^{\infty}} + \gamma T_1 \|u_0\|_{L^2} \right)^{\frac{1}{2}} = \log\left( 1 + \frac{2}{\varepsilon} \right)
\end{equation}
and assume that there is a $x_{0} \in \mathbb{S}$ such that
\begin{equation}
\label{derivative-bound}
u'_{0}(x_{0}) \leq - (1 + \epsilon) \sqrt {\gamma}
\left ( \|u_0\|_{L^{\infty}} + \gamma T_1 \| u_0 \|_{L^2} \right )^{\frac {1}{2}}.
\end{equation}
Then the solution $u(t,x)$ in Lemma \ref{lemma-wellposedness} blows
up in a finite time $T \in (0,T_1)$ in the sense of Lemma \ref{lemma-blowup}.
\end{theorem}

\begin{proof}
Define $V(t,\xi) = u_x(t,X(t,\xi))$. By Lemmas
\ref{lemma-wellposedness} and \ref{lemma-characteristics},
$V(t,\xi)$ is absolutely continuous on $[0,T) \times \mathbb{S}$
and almost everywhere differentiable on $(0,T) \times \mathbb{S}$,
so that
\begin{eqnarray*}
\dot{V} = \left( u_{tx} + u u_{xx} \right) \biggr|_{x = X(t,\xi)}
= \left( \gamma u - u_x^2 \right) \biggr|_{x = X(t,\xi)} = -V^2 +
\gamma U, \quad {\rm a.e.} \quad \xi \in \mathbb{S}, \quad t \in
(0,T).
\end{eqnarray*}
By Lemma \ref{lemma-bounds}, we obtain the apriori differential estimate
\begin{equation}
\label{diff-inequality}
\dot{V} \leq - V^2 + \gamma \left( \| u_0 \|_{L^{\infty}} + \gamma t \| u_0 \|_{L^2} \right),
 \quad {\rm a.e.} \quad \xi \in \mathbb{S}, \quad t \in
(0,T).
\end{equation}
Since $u_0'(x)$ is a continuous, mean-zero, periodic function of $x$ on $\mathbb{S}$
and assumption (\ref{derivative-bound}) is satisfied for fixed $\varepsilon > 0$,
there exists $\tilde{x}_0$ such that
$$
V(0,\tilde{x}_0) = -(1+\varepsilon) h(T_{1}),
$$
where
$$
h(T_{1})= \sqrt {\gamma} \left ( \|u_0\|_{L^{\infty}} + \gamma T_1 \|u_0\|_{L^2}\right )^{\frac{1}{2}}.
$$
Thanks to the apriori estimate (\ref{diff-inequality}), $V(t) := V(t,\tilde{x}_0)$ satisfies
\begin{equation}
\label{inequality} \left\{ \begin{array}{l} \dot{V}(t) \leq
-V^2(t) + h^{2}(T_{1}), \quad {\rm a.e.} \quad t\in [0,T_{1}]\cap (0,T), \\
V(0) = -(1+\varepsilon) h(T_1). \end{array} \right.
\end{equation}
By the comparison principle for ODEs, we have
$$
V(t) \leq V_+(t) < 0, \quad t\in [0,T_{1}]\cap[0,T),
$$
where $V_+(t)$ solves the equation
\begin{equation}
\label{super-solution-V} \left\{ \begin{array}{l} \dot{V}_+(t) =
-V_+^2(t) + h^{2}(T_{1}), \quad   t \in [0,T_{1}), \\ V_+(0) =
V(0). \end{array} \right.
\end{equation}
Equation (\ref{super-solution-V}) admits an implicit solution
$$
\frac{V_+(t)+h(T_{1})}{V_+(t)-h(T_{1})} = \frac{V(0)+h(T_{1})}{V(0)-h(T_{1})} e^{2h(T_{1})t}, \quad t \in [0,T_1).
$$
If $T_1$ is the smallest positive root of (\ref{time-bound}), then
$$
\frac{V_+(t)+h(T_{1})}{V_+(t)-h(T_{1})} = \frac{\varepsilon}{2 + \varepsilon} e^{2h(T_{1}) t}
\uparrow 1, \quad \mbox{\rm as} \quad t \uparrow T_1,
$$
so that $\lim_{t \uparrow T_1} V_+(t) = -\infty$. Therefore, there is $T \in (0,T_1)$ such that  $\displaystyle\lim_{t\uparrow T} V(t)=-\infty$.
\end{proof}

\begin{remark}
Note that if $\varepsilon \to \infty$ and the assumption
of Theorem \ref{theorem-blowup-characteristics} still holds, then $T \to 0$. This means that the steeper the
slope of the initial data $u_0(x)$ is, the quicker the solution $u(t,x)$ blows up.
\end{remark}

\begin{remark}
Since $\| u(t,\cdot) \|_{L^{\infty}}$ remains bounded on $[0,T)$ thanks to Lemma \ref{lemma-bounds},
the blow-up of Theorem \ref{theorem-blowup-characteristics} corresponds to criterion (\ref{wave-breaking-criterion}) of the wave breaking in the Ostrovsky--Hunter equation (\ref{OstHunter}).
\end{remark}

By Theorem \ref{theorem-blowup-characteristics}, we have the following two corollaries.

\begin{corollary}
\label{corollary-1}
Assume that $u_0(x) \in H^s(\mathbb{S})$, $s \ge 2$ is even and
non-constant. Then for sufficiently large $n$, the corresponding
solution $u(t,x)$ to the Cauchy problem (\ref{cauchy})
with initial data $u_0(nx)$ blows up in finite time.
\end{corollary}

\begin{proof}
Take $x_0\in \mathbb{S}$ such that $u'_0(x_0) = \inf_{x\in \mathbb{S}} u'_0(x)$. Since
$u_0 \in C^1(\mathbb{S})$ is even and periodic, it follows that
$u'_0(x_0)  \le 0$ and $\sup_{x\in \mathbb{S}} u'_0(x) = -u_0'(x_0) \geq 0$.
Thus, we deduce that
\begin{equation*}
\left( \inf_{x\in \mathbb{S}}u'_0(x)\right)^2=\left( \sup_{x\in
\mathbb{S}}u'_0(x)\right)^2>\int_{\mathbb{S}}(u'_0(x))^2dx.
\end{equation*}
Let $\tilde{u}_0(x) = u_0(nx)$ for a positive integer $n$.
Thanks to $1$-periodicity of $u_0(x)$, we have $\| \tilde{u}'_0 \|_{L^2} = n \|u'_0\|_{L^2}$,
$\| \tilde{u}_0 \|_{L^2}=\|u_0\|_{L^2}$, and
$\| \tilde{u}_0 \|_{L^{\infty}} = \|u_0\|_{L^{\infty}}$. From the above
inequality, we see that the assumption of Theorem \ref{theorem-blowup-characteristics}
for any $\varepsilon > 0$ is satisfied by the initial data
$\tilde{u}_0(x) = u_0(nx)$ provided $n$ is large enough.
\end{proof}

\begin{corollary}
\label{corollary-2}
Assume that $u_0(x) \in H^s(\mathbb{S})$, $s \ge 2$ and
$|\inf_{x\in\mathbb{S}}u'_0(x)|\geq |\sup_{x\in\mathbb{S}}u'_0(x)| > 0$.
Then for sufficiently large $n$, the corresponding solution $u(t,x)$
to the Cauchy problem (\ref{cauchy})  with initial data $u_0(nx)$ blows up in finite time.
\end{corollary}

\begin{proof}
The assumption and the mean value theorem imply that there is a point
$x_0\in \mathbb{S}$ such that
$$
\left( \inf_{x\in \mathbb{S}} u'_0(x)\right)^2 \geq
\left(\sup_{x\in\mathbb{S}}u'_0(x)\right)^2 > (u'_0(x_0))^2
\geq \int_{\mathbb{S}}(u'_0(x))^2dx.
$$
Thus, we can obtain the desired result in view of the proof of
Corollary \ref{corollary-1}.
\end{proof}

Our final result specifies the rate at which the wave breaks
in the Cauchy problem (\ref{cauchy}). We use again the fact
that the blow-up time $T$ is independent of $s \geq 2$ for the solution $u(t,x)$
in Lemma \ref{lemma-wellposedness}, so that the initial data $u_0(x)$
can be considered in $H^3(\mathbb{S})$.

\begin{theorem}
\label{theorem-blowup-rate}
Let $u_0(x) \in H^3(\mathbb{S})$ and $T \in (0,\infty)$ be
the finite blow-up time of the solution $u(t,x)$ in Lemma \ref{lemma-wellposedness}.
Then we have
\begin{equation}
\label{result-1}
\lim_{t \uparrow T} \left( (T-t)
\inf_{x\in\mathbb{S}} u_{x}(t,x) \right)=-1
\end{equation}
and
\begin{equation}
\label{result-2}
\lim_{t \uparrow T}\left( (T-t)
\sup_{x\in\mathbb{S}} u_{x}(t,x) \right)= 0.
\end{equation}
\end{theorem}

\begin{proof}
Let $m(t) := \inf_{x\in\mathbb{S}} u_{x}(t,x)$. By the assumption of $T<\infty$ and
Lemma \ref{lemma-blowup}, we have $\lim_{t \uparrow T} m(t)=-\infty$.
By Theorem 2.1 in Constantin \& Escher \cite{C-E}, for every $t\in [0,T)$, there
exists at least one point $\xi(t)\in \mathbb{S}$ such that
$m(t):= u_{x}(t,\xi(t))$ and $u_{xx}(t,\xi(t))=0$. Moreover, $m(t)$ (and $\xi(t)$) is
absolutely continuous on $[0,T)$, almost everywhere differentiable on $(0,T)$, and satisfies
\begin{equation}
\label{equality-again}
\frac{d}{dt} m(t) = u_{tx}(t,\xi(t)) = -m^{2}(t)+ \gamma u(t, \xi(t)) \quad
{\rm a.e.} \quad t \in (0,T).
\end{equation}
Set $K(T)=  \gamma \left ( \| u_{0}\|_{L^{\infty}} + \gamma T \| u_{0}\|_{L^{2}} \right)$.
By Lemma \ref{lemma-bounds}, we obtain
\begin{equation}
\label{inequality-again} -m^{2}(t)-K(T) \leq \frac{d}{dt}m(t)\leq
-m^{2}(t)+K(T)\quad {\rm a.e.} \quad t\in (0,T).
\end{equation}
Let us now choose $\varepsilon\in (0,1)$. Since $\lim_{t \uparrow T} m(t)=-\infty$, one can find
$t_{0}\in [0, T)$ such that
$$
m(t_0)<-\sqrt{K(T)+\frac{K(T)}{\varepsilon}}.
$$
By the continuation of solutions of (\ref{inequality-again})
and the absolute continuity of $m(t)$, it follows
that $m$ is decreasing on $[t_0,T)$ so that
$$
m(t) \leq m(t_0) < -\sqrt{K(T)+\frac{K(T)}{\varepsilon}}<-\sqrt{\frac{K(T)}{\varepsilon}},\quad t\in[t_0,T)
$$
and
\begin{equation*}
1-\varepsilon\leq \frac{d}{dt}\left(\frac{1}{m(t)}\right)\leq
1+\varepsilon.
\end{equation*}
Integrating the above relation on $(t,T)$ with $t\in [t_0,T)$ and
noticing that $\lim_{t\uparrow T}m(t)=-\infty$, we deduce that
\begin{equation*}
(1-\varepsilon)(T-t)\leq -\frac{1}{m(t)}\leq (1+\varepsilon)(T-t).
\end{equation*}
Since $\varepsilon\in(0,1)$ is arbitrary, in view of the definition
of $m(t)$, the above inequality in the limit $\varepsilon \downarrow 0$
implies the desired result (\ref{result-1}).

Now let $M(t) := \sup_{x \in \mathbb{S}} u_x(t,x)$.
By the same Theorem 2.1 in Constantin \& Escher \cite{C-E},
for every $t\in [0,T)$, there
exists at least one point $\eta(t)\in \mathbb{S}$ such that
$M(t)  = u_x(t,\eta(t))$ and $u_{xx}(t,\eta(t))=0$.
Repeating the same arguments, we have
$$
\frac{d}{dt} M(t) =  -M^2(t) + \gamma u (t, \eta(t)) \le \gamma
\left ( \|u_0 \|_{L^{\infty}} + \gamma t \|u_0\|_{L^2} \right), \quad
\mbox{for all} \quad t\in (0,T),
$$
so that
\begin{equation}
\label{bound-final}
M(t) \le \sup_{x \in \mathbb{S}}u_0'(x) + \gamma \left( T \|u_0 \|_{L^{\infty}}
+ \frac{\gamma T^2}{2} \|u_0 \|_{L^2} \right) < + \infty.
\end{equation}
Since $u(t,x)$ is periodic on $\mathbb{S}$ for all $t \in [0,T)$ and belong
to $C([0,T);H^3(\mathbb{S}))$,
there exists $\xi_0(t) \in \mathbb{S}$ for every $t \in [0, T)$ such that
$u_x(t, \xi_0(t)) = 0 $. Therefore, $M(t) \ge u_x(t, \xi_0(t)) = 0$ for all $t \in [0,T)$,
so that bound (\ref{bound-final})
yields the desired result (\ref{result-2}). This completes the proof of the theorem.
\end{proof}

\section{Wave breaking on an infinite line}

To extend our results on wave breaking in the Ostrovsky--Hunter equation (\ref{OstHunter}) from
a circle $\mathbb{S}$ to an infinite line $\mathbb{R}$, we are going to use an additional conserved
quantity of the Ostrovsky--Hunter equation. Consider the Cauchy problem in the form
\begin{equation}
\label{cauchy-infinite}
\left\{ \begin{array}{l} u_t + u u_x = \gamma \partial_x^{-1} u, \quad  x \in \mathbb{R}, \;\;
t > 0, \\
u(0,x) = u_0(x), \qquad \phantom{t} x \in \mathbb{R}, \end{array} \right.
\end{equation}
where $\gamma > 0$ and $\partial_x^{-1} u = \int_{-\infty}^x u(t,x') dx'$. To control
$\partial_x^{-1} u$, we define $\| u \|_{\dot{H}^{-1}} := \| \partial_x^{-1} u \|_{L^2}$. The
local well-posedness result is given by the following lemma.

\begin{lemma}  Assume that $u_{0}(x) \in H^{s}(\mathbb{R}) \cap \dot{H}^{-1}(\mathbb{R})$, $s \ge 2$.
Then  there exist a maximal time $T=T(u_{0})>0$  and a unique solution $u(t,x)$ to the
Cauchy problem (\ref{cauchy-infinite}) such that
$$
u(t,x) \in C([0,T);H^{s}(\mathbb{R})\cap \dot{H}^{-1}(\mathbb{R}))\cap
C^{1}([0,T);H^{s-1}(\mathbb{R}))
$$
with the following three conserved quantities
\begin{equation}
\label{mass-conservation}
\int_{\mathbb{R}} u(t, x) dx  = 0,  \quad t \in [0, T),
\end{equation}
\begin{equation}
\label{power-conservation}
Q = \int_{\mathbb{R}} u^2(t, x) dx = \int_{\mathbb{R}} u_0^2(x) dx,    \quad  t \in [0, T),
\end{equation}
and
\begin{equation}
\label{energy-conservation}
E = \int_{\mathbb{R}} \left[ \gamma (\partial_x^{-1} u)^2 + \frac{1}{3} u^3 \right] dx =
\int_{\mathbb{R}} \left[ \gamma (\partial_x^{-1} u_0)^2 + \frac{1}{3} u_0^3 \right] dx, \quad t \in [0, T).
\end{equation}
Moreover, the solution depends continuously on the initial data,
i.e. the mapping $u_{0}\mapsto u : \; H^{s}(\mathbb{R})
\rightarrow C([0,T);H^{s}(\mathbb{R}) \cap \dot{H}^{-1}(\mathbb{R}))\cap
C^{1}([0,T);H^{s-1}(\mathbb{R}))$ is continuous.
\label{lemma-wellposedness-infinite}
\end{lemma}

\begin{proof}
If $u_0(x) \in H^s(\mathbb{R}) \cap \dot{H}^{-1}(\mathbb{R})$, then $\partial_x^{-1} u_0(x) \in H^3(\mathbb{R})$,
so that $\int_{\mathbb{R}} u_0(x) dx = 0$. By the main theorem in \cite{ScWa} in the context
of the short-pulse equation (\ref{short-pulse}), existence, uniqueness and continuous dependence
of the solution $u(t,x) \in C([0,T);H^s(\mathbb{R})) \cap C^{1}([0,T);H^{s-1}(\mathbb{R}))$ is proved,
so that
$$
\gamma \partial_x^{-1} u(t,x) = u_t(t,x) + u(t,x) u_x(t,x) \in C((0,T);H^{s-1}(\mathbb{R})).
$$
Therefore, $u(t,x) \in C([0,T);H^s(\mathbb{R})\cap \dot{H}^{-1}(\mathbb{R}))$ in view of continuity
of $\| u \|_{\dot{H}^{-1}}$ as $t \downarrow 0$.
Because $f \in H^1(\mathbb{R})$ implies $\lim_{|x| \to \infty} f(x) = 0$, the zero-mass constraint (\ref{mass-conservation}) holds. Let us define
$$
\gamma \partial_x^{-2} u(t,x) = \left(\partial_x^{-1} u(t,x) \right)_t + \frac{1}{2} u^2(t,x).
$$
By uniqueness of the strong solution $u(t,x)$ satisfying the constraint (\ref{mass-conservation}) on $[0,T)$,
we obtain
$$
\lim_{|x| \to \infty} \partial_x^{-2} u(t,x) = 0, \quad \forall t \in [0,T).
$$
Using balance equations for the densities of $Q$ and $E$, we write
\begin{eqnarray*}
(u^2)_t &  = & \left( \gamma (\partial_x^{-1} u)^2 - \frac{2}{3} u^3
\right)_x, \quad x \in \mathbb{R}, \; t \in (0,T), \\
\left[  \gamma (\partial_x^{-1} u)^2 + \frac{1}{3} u^3 \right]_t & = & \left[
\gamma^2 (\partial_x^{-2} u)^2 - \frac{1}{4} u^4 \right]_x, \quad x \in \mathbb{R}, \; t \in (0,T).
\end{eqnarray*}
Integrating the balance equation in $x$ on $\mathbb{R}$ for any $t \in (0,T)$, we obtain
conservation of $Q$ and $E$. Their initial values as $t \downarrow 0$
are computed from the initial condition $u_0(x)$ thanks to the fact that
$u_0(x) \in H^s(\mathbb{R}) \cap \dot{H}^{-1}(\mathbb{R})$.
\end{proof}

The blow-up alternative in Lemma \ref{lemma-blowup} holds on an
infinite line thanks to Sobolev's embedding $H^3(\mathbb{R})$
into $C^2(\mathbb{R})$ and the density arguments. Since the
application (\ref{Holder-inequality}) of the H\"{o}lder inequality
is not valid on $\mathbb{R}$, Theorem \ref{theorem-wave-breaking}
can not be extended on an infinite line. However, we can still use
the method of characteristics and extend Theorems
\ref{theorem-blowup-characteristics} and \ref{theorem-blowup-rate}
from $\mathbb{S}$ to $\mathbb{R}$.

For all $t \in [0,T)$ and $\xi \in \mathbb{R}$, we define
$$
x = X(t,\xi), \quad u(t,x) = U(t,\xi), \quad \partial_x^{-1} u(t,x) = G(t,\xi),
$$
so that the same system (\ref{characteristics}) is considered. Lemma \ref{lemma-characteristics}
holds on $\mathbb{R}$, while Lemma \ref{lemma-bounds} is replaced with the following lemma.

\begin{lemma}
Let $u_{0}(x) \in H^s(\mathbb{R}) \cap \dot{H}^{-1}(\mathbb{R})$, $s \ge 2$ and $T>0$ be the
maximal existence time of the solution $u(t,x)$ in Lemma \ref{lemma-wellposedness-infinite}.
The solution $u(t,x)$ satisfies
$$
\sup_{s \in [0,t]} \| u(s, \cdot)\|_{L^{\infty}} \leq \| u_{0}\|_{L^{\infty}} + C t + \frac{\gamma Q}{6} t^2,
\quad t\in [0,T),
$$
where
\begin{eqnarray*}
C = \frac{\sqrt{\gamma}}{\sqrt{2}} \left( E + \gamma Q + \frac{1}{3} Q \| u_0 \|_{L^{\infty}} \right)^{\frac{1}{2}}.
\end{eqnarray*}
\label{lemma-bounds-infinite}
\end{lemma}

\begin{proof}
From conserved quantities (\ref{power-conservation}) and (\ref{energy-conservation}), we obtain
\begin{eqnarray*}
\| \partial_x^{-1} u(t,\cdot) \|^2_{H^1} = \| u(t,\cdot) \|^2_{L^2} + \| \partial_x^{-1} u(t,\cdot) \|^2_{L^2} & = & Q +
\frac{1}{\gamma} \left( E - \frac{1}{3} \int_{\mathbb{R}} u^3(t,x) dx \right) \\
& \leq & Q +
\frac{1}{3 \gamma} \left( 3 E + Q \| u(t,\cdot) \|_{L^{\infty}} \right).
\end{eqnarray*}
Let $S(t) := \sup_{s \in [0,t]} \| u(s, \cdot)\|_{L^{\infty}}$.
Thanks to Sobolev's embedding of $H^1$ into $L^{\infty}$, we have
\begin{eqnarray*}
S(t) & \leq & \| u_0 \|_{L^{\infty}} + \gamma t \sup_{s \in [0,t]}
\| \partial_x^{-1} u(s,\cdot) \|_{L^{\infty}} \\
& \leq & \| u_0 \|_{L^{\infty}} + \gamma t \left(
\frac{Q}{2} +
\frac{1}{6 \gamma} \left( 3 E + Q S(t) \right) \right)^{\frac{1}{2}},
\end{eqnarray*}
from which the bound on $S(t)$ is proved after algebraic manipulations.
\end{proof}

Theorem \ref{theorem-blowup-characteristics} is extended to the infinite line
in the following theorem.

\begin{theorem}
\label{theorem-blowup-characteristics-infinite}
Let $\varepsilon>0$ and $u_{0}(x) \in H^{s}(\mathbb{R}) \cap \dot{H}^{-1}(\mathbb{R})$, $s \ge 2$.
Let $T_1$ be the smallest positive root of the equation
$$
2 \sqrt{\gamma} T_1 \left(\| u_{0}\|_{L^{\infty}} + C T_1 +
\frac{\gamma Q}{6} T_1^2 \right)^{\frac{1}{2}}
= \log\left( 1 + \frac{2}{\varepsilon} \right)
$$
and assume that there is a $x_{0} \in \mathbb{S}$ such that
$$
u'_{0}(x_{0}) \leq - (1 + \epsilon) \sqrt {\gamma}
\left ( \| u_{0}\|_{L^{\infty}} + C T_1 + \frac{\gamma Q}{6} T_1^2 \right )^{\frac {1}{2}},
$$
where $C$ is defined in Lemma \ref{lemma-bounds-infinite}.
Then the solution $u(t,x)$ in Lemma \ref{lemma-wellposedness-infinite} blows
up in a finite time $T \in (0,T_1)$.
\end{theorem}

\begin{proof}
The proof is similar to that of Theorem \ref{theorem-blowup-characteristics}.
\end{proof}

Finally, Theorem \ref{theorem-blowup-rate} remains valid on an infinite line since the
proof does not depend on the definition of $K(T)$.

\section{Numerical evidence of wave breaking}

\begin{figure}
\begin{center}
\includegraphics[width=0.7\textwidth]{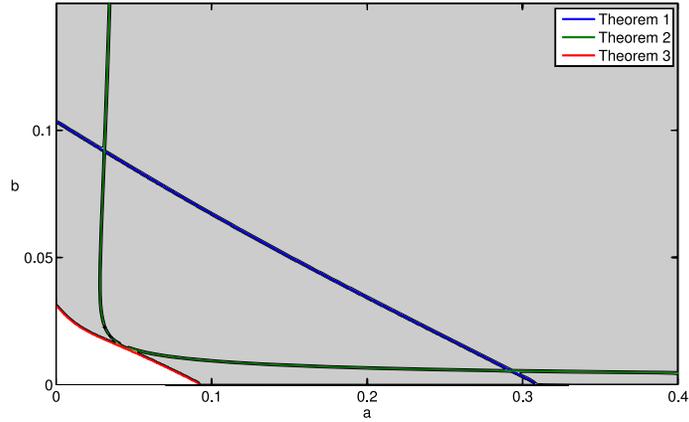}
\end{center}
\caption{Lower bounds of the domains for the wave breaking
conditions of Theorems \ref{theorem-hunter},
\ref{theorem-wave-breaking}, and
\ref{theorem-blowup-characteristics}. Shaded area shows where the
wave breaking condition of Theorem
\ref{theorem-blowup-characteristics} is satisfied.}
\label{fig-parameter-plane}
\end{figure}

We consider the Cauchy problem (\ref{cauchy}) on
$\mathbb{S}$ for $\gamma = 1$ and the initial data
\begin{equation} \label{initial-data}
u_0(x) = a \cos(2\pi x) + b \sin(4\pi x),
\end{equation}
where $(a,b)$ are parameters. Using elementary calculus, we compute
$$
\int_{\mathbb{S}} (u_0'(x))^3 dx = - 12 \pi^3 a^2 b, \quad
\int_{\mathbb{S}} u_0^2(x) dx = \frac{1}{2} (a^2+b^2)
$$
and
\begin{eqnarray*}
m & = & -\inf_{x \in \mathbb{S}} u_0'(x) = 2 \pi (a+2b), \\
M & = & \sup_{x \in \mathbb{S}} |u_0(x)| = \frac{1}{4}
\left(3a + \sqrt{a^2+32b^2} \right)
\sqrt{1 - \left( \frac{-a+\sqrt{a^2+32b^2}}{8b} \right)^2}.
\end{eqnarray*}
Figure \ref{fig-parameter-plane} compares the theoretical
estimates of the wave breaking regions on the quarter-plane $(a,b)
\in \mathbb{R}_+^2$. According to Theorem \ref{theorem-hunter},
the blow-up occurs under the condition $m^3 > 4M (4 + m)$. The
lower bound of the domain, where $m^3 > 4M (4 + m)$, is shown by
the blue line on Figure \ref{fig-parameter-plane}. According to
Theorem \ref{theorem-wave-breaking}, two criteria
(\ref{condition-one}) and (\ref{condition-two}) exist. The lower
bound of the domain, where the first criterion
(\ref{condition-one}) is met, is shown on Figure
\ref{fig-parameter-plane} by the green curve. The domain of the
second criterion (\ref{condition-two}) is, however, beyond the
scale of Figure \ref{fig-parameter-plane}. Indeed, the latter
domain corresponds to the quarter circle on $(a,b)$-plane with the
radius $\frac{3}{2\sqrt{2}} \approx 1.06$. Finally, the lower
bound of the domain given by the criterion of Theorem
\ref{theorem-blowup-characteristics} is shown on Figure
\ref{fig-parameter-plane} by the red line. We can see from the
figure that the latter criterion (\ref{condition-one}) is the
sharpest one with the largest wave breaking region shown by shaded
area on Figure \ref{fig-parameter-plane}.

\begin{figure}
\begin{center}
\begin{tabular}{cc}
\includegraphics[width=0.45\textwidth]{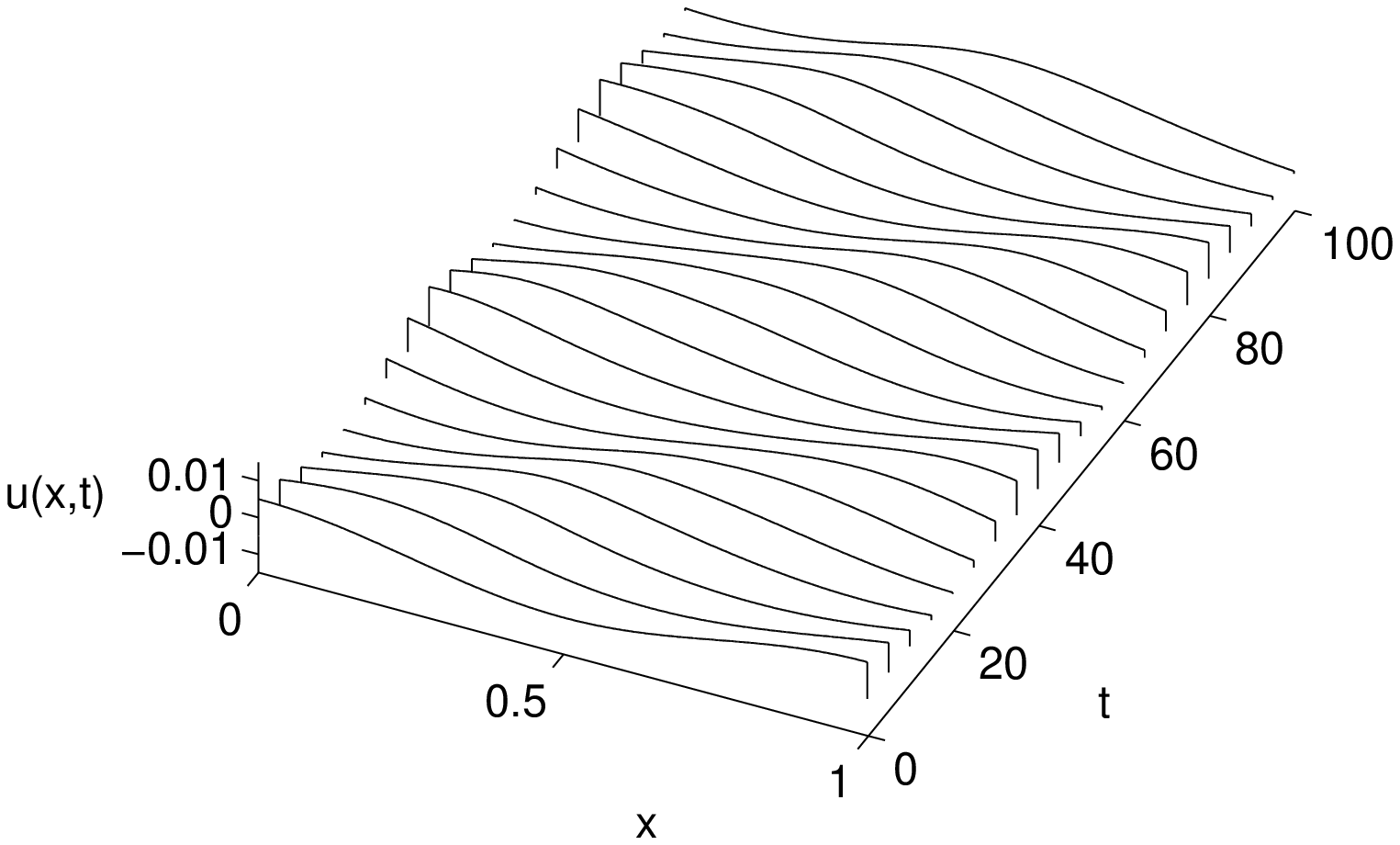} & \includegraphics[width=0.45\textwidth]{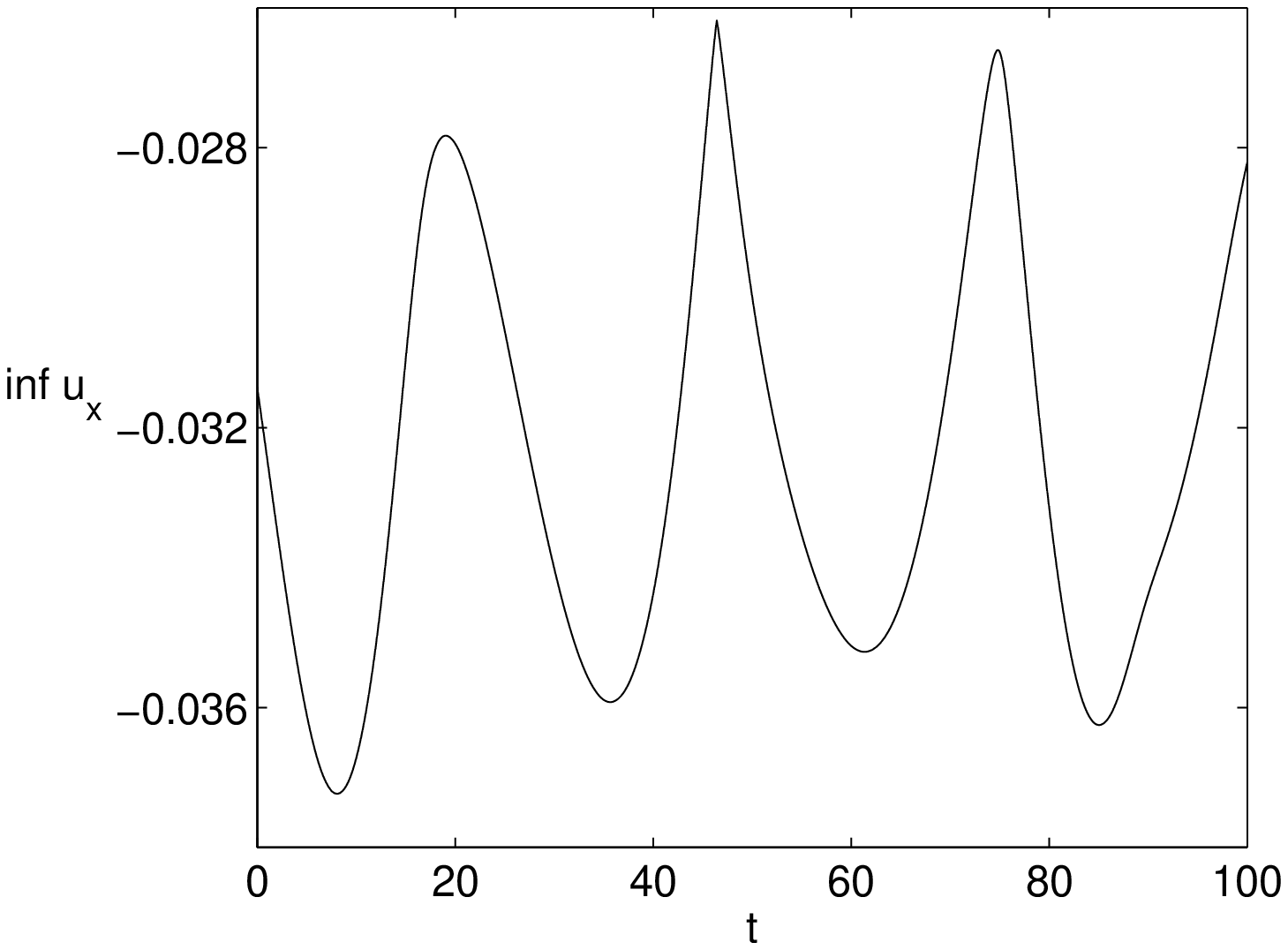}\\
\includegraphics[width=0.45\textwidth]{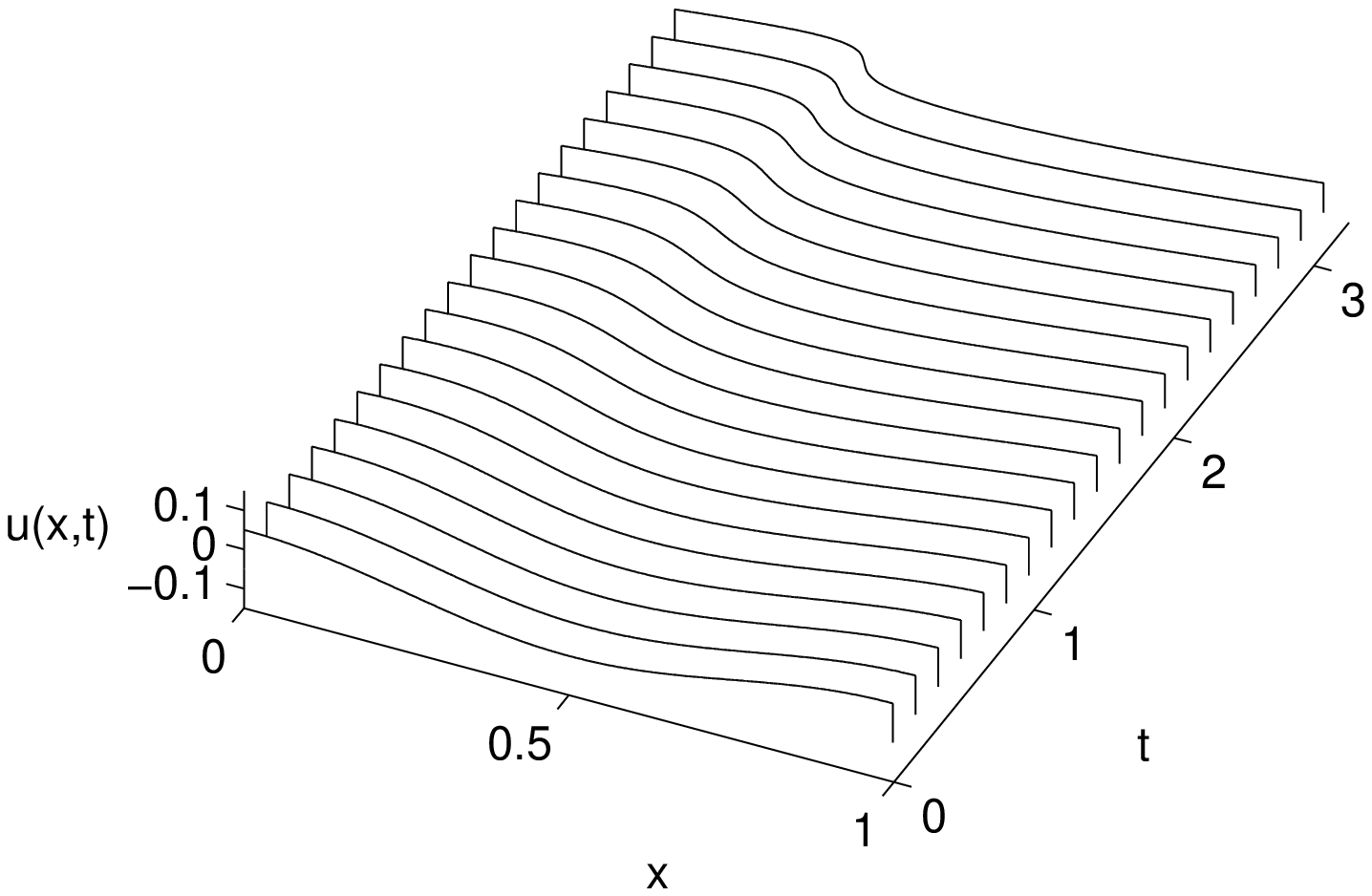} & \includegraphics[width=0.45\textwidth]{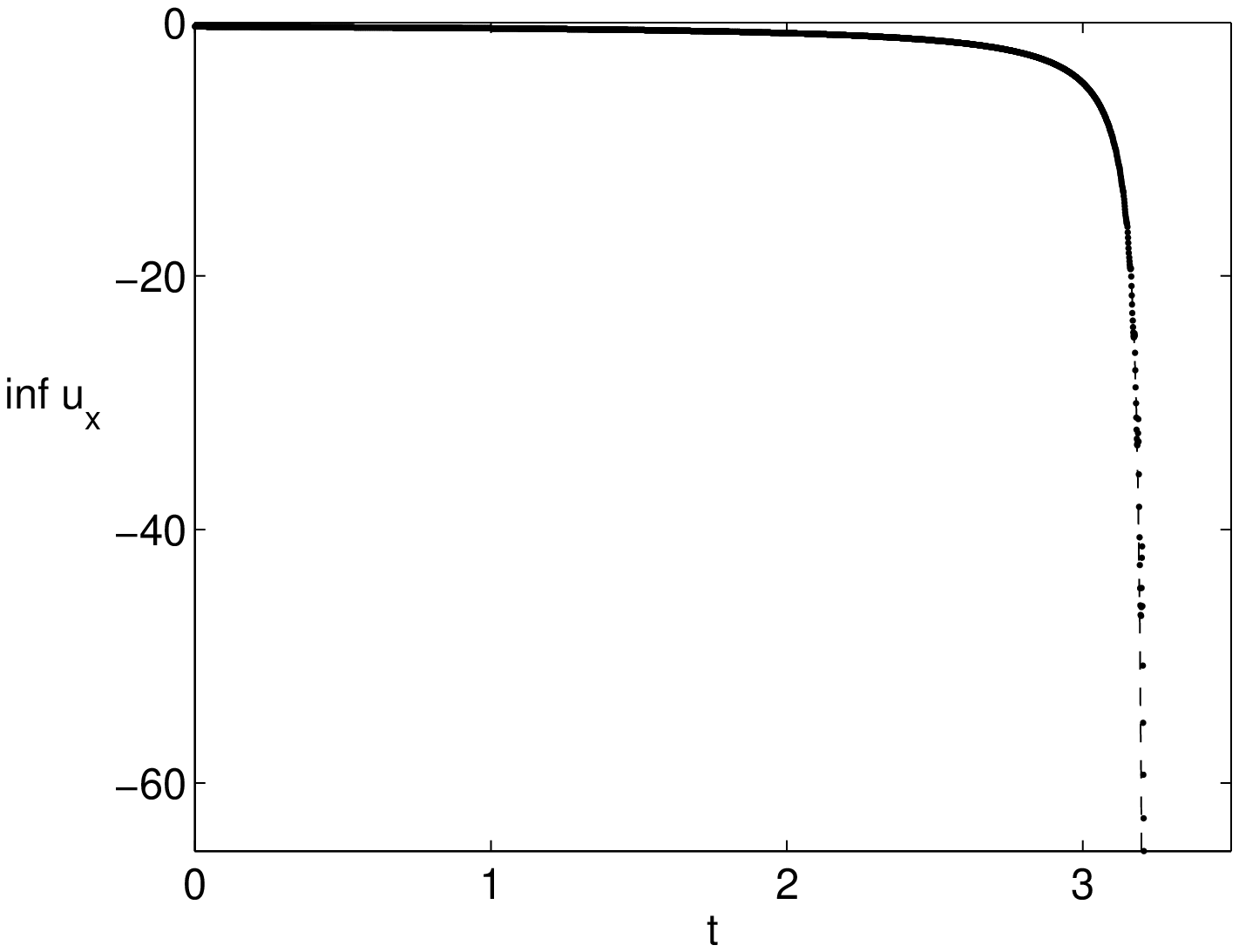}\\
\end{tabular}
\end{center}
\caption{Solution surface $u(t,x)$ (left) and $\inf_{x \in
\mathbb{S}} u_x(t,x)$ versus $t$ (right) for two simulations with
$a=0.005$, $b=0$ (top) and $a=0.05$, $b=0$ (bottom). The dashed
curve on the bottom right panel shows the least squares fit
$-1/(B+Ct)$ with $C \approx -1.009$ and $B \approx 3.213$.}
\label{fig-1}
\end{figure}

\begin{figure}
\begin{center}
\begin{tabular}{cc}
\includegraphics[width=0.45\textwidth]{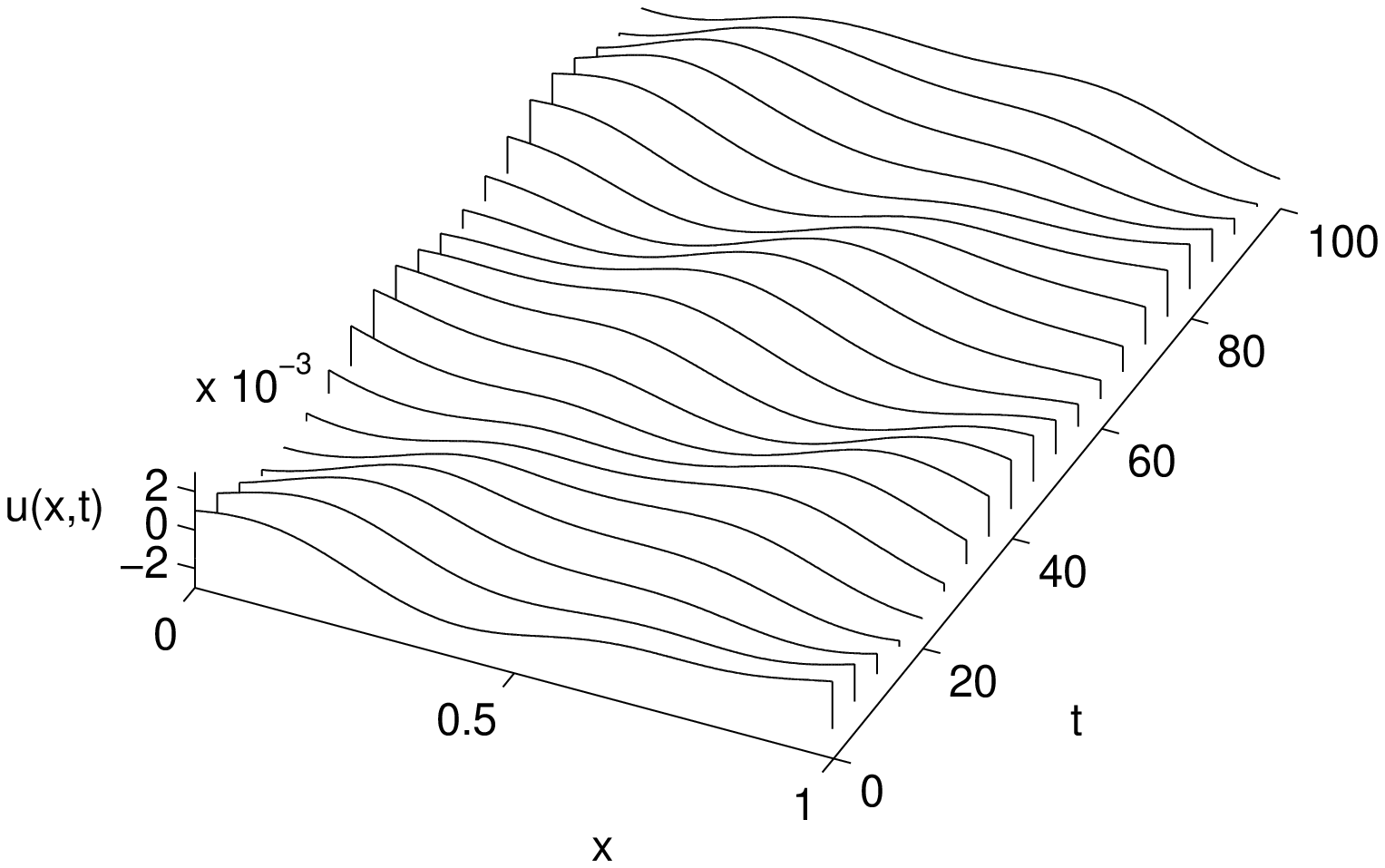} & \includegraphics[width=0.45\textwidth]{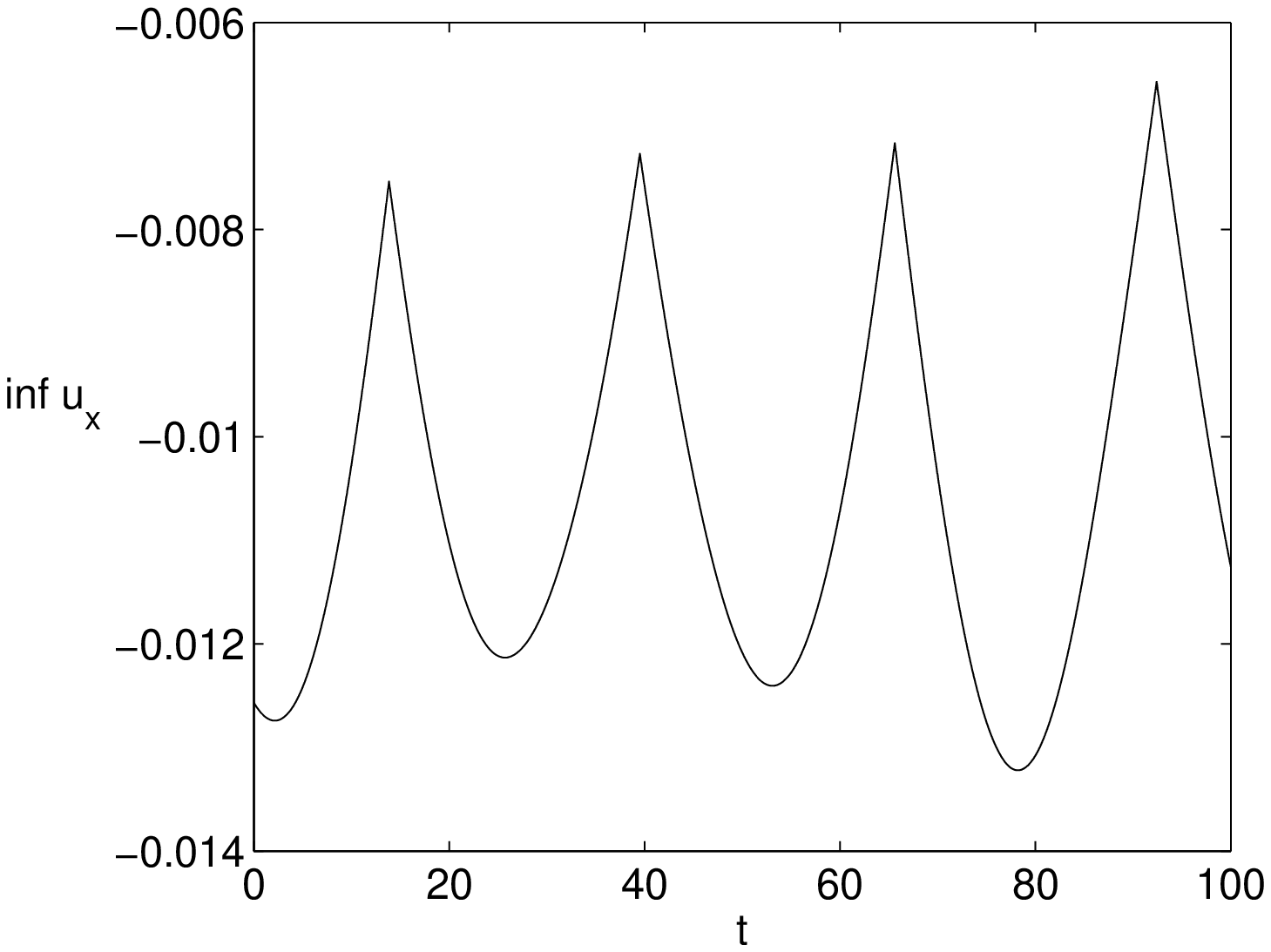}\\
\includegraphics[width=0.45\textwidth]{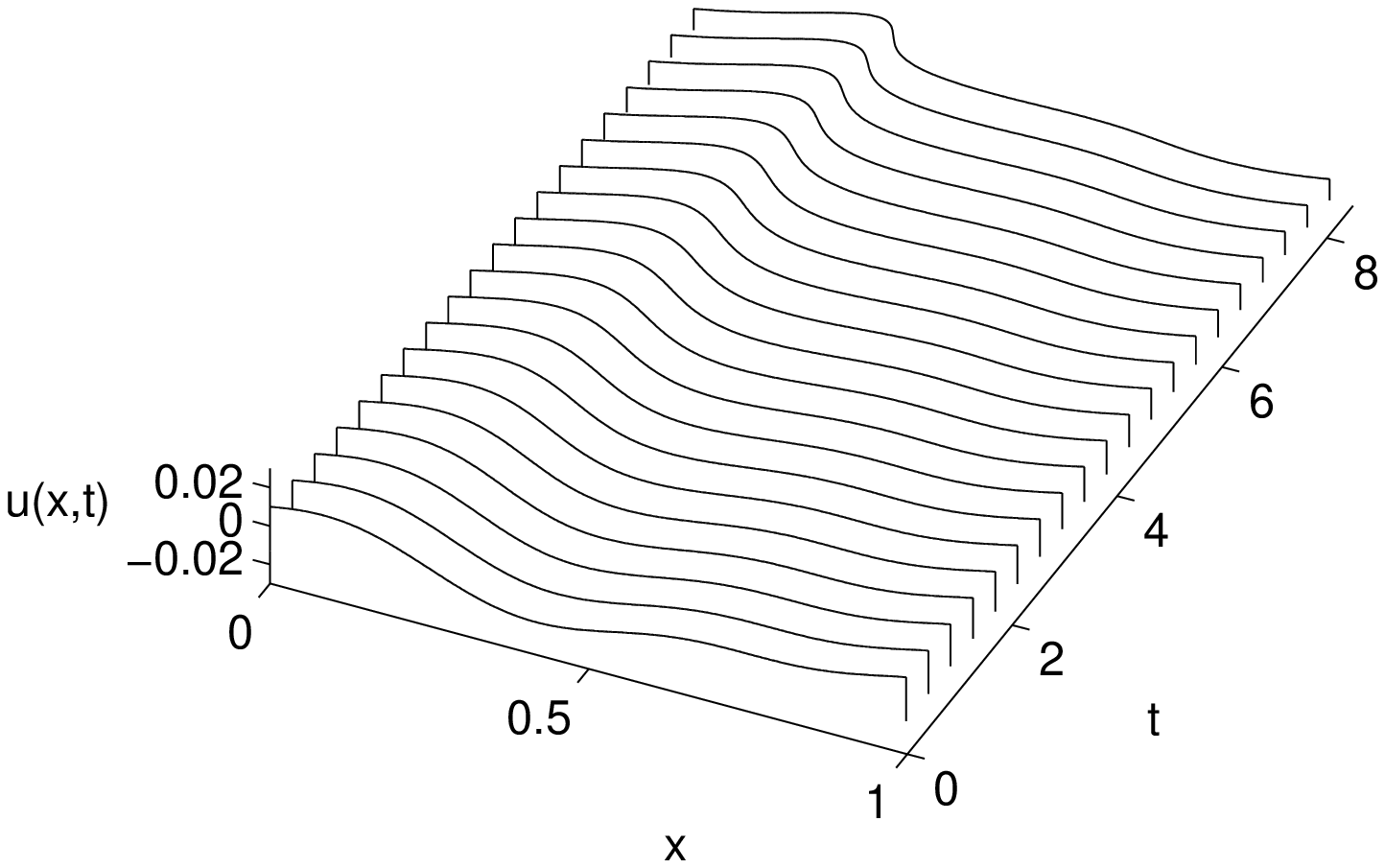} & \includegraphics[width=0.45\textwidth]{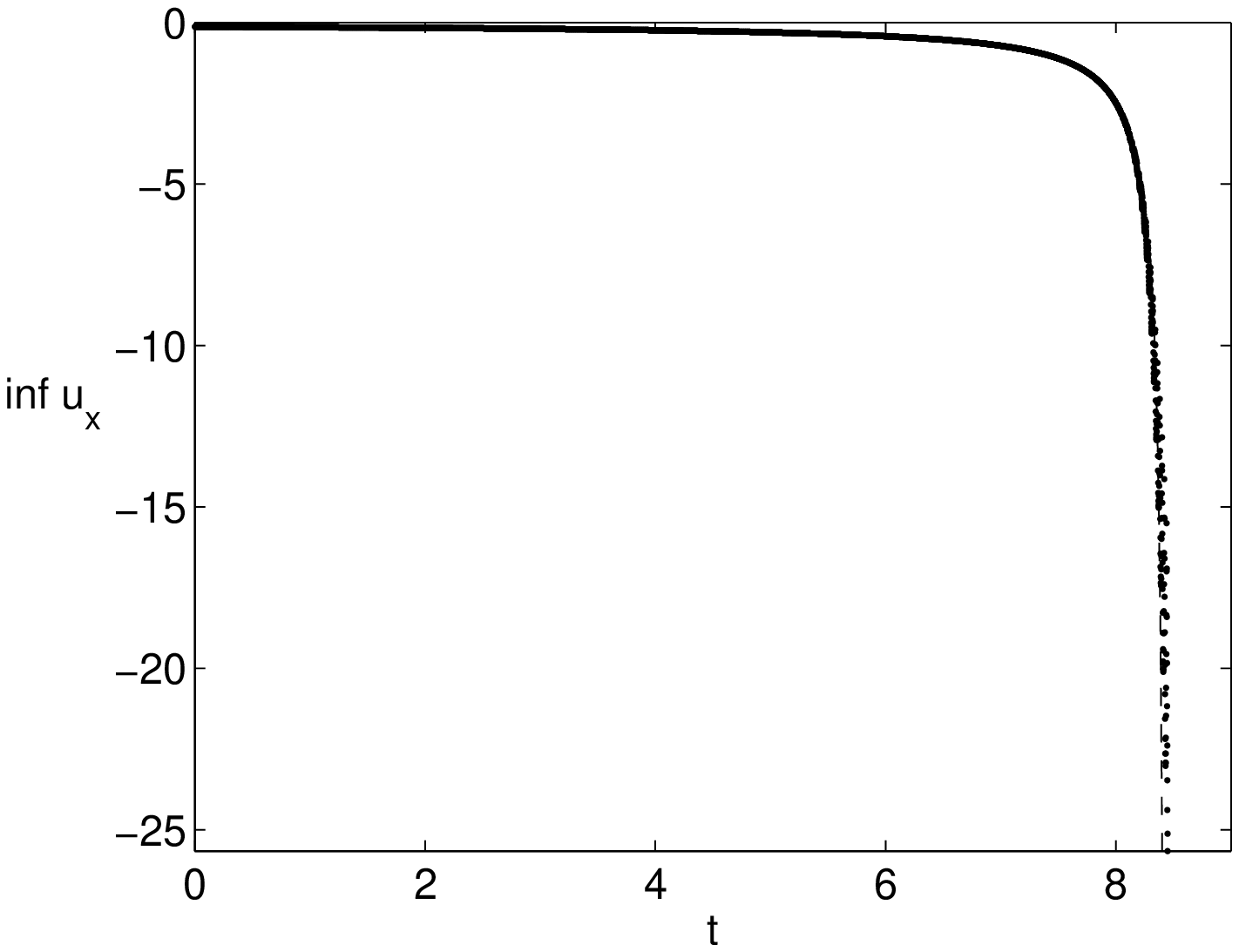}\\
\end{tabular}
\end{center}
\caption{The same as Figure \ref{fig-1} but for $a=0.001$,
$b=0.0005$ (top) and $a=0.01$, $b=0.005$ (bottom). The least
squares fit is computed with $C \approx -1.042$ and $B \approx
8.442$.} \label{fig-2}
\end{figure}

\begin{figure}
\begin{center}
\begin{tabular}{cc}
\includegraphics[width=0.45\textwidth]{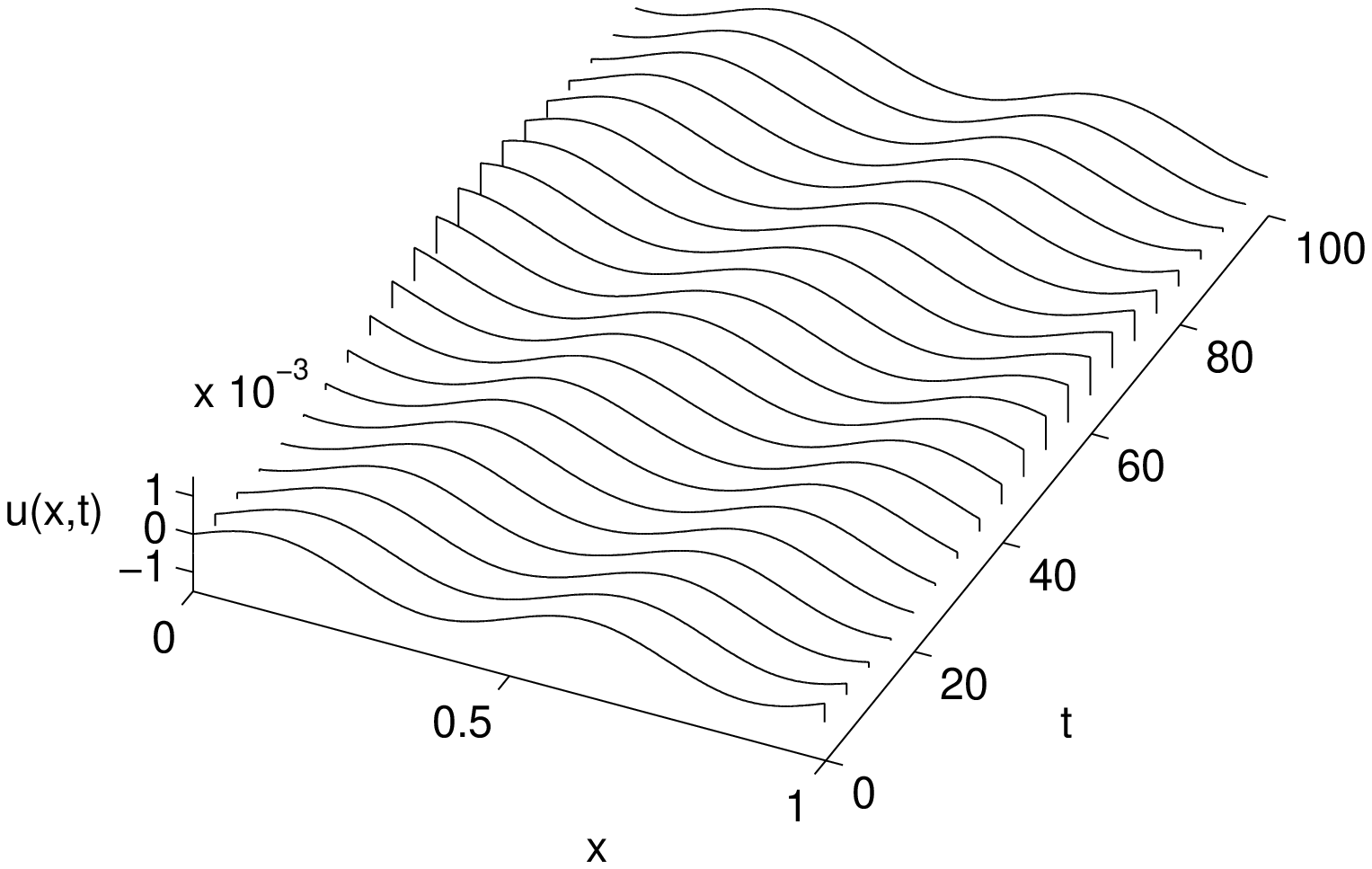} & \includegraphics[width=0.45\textwidth]{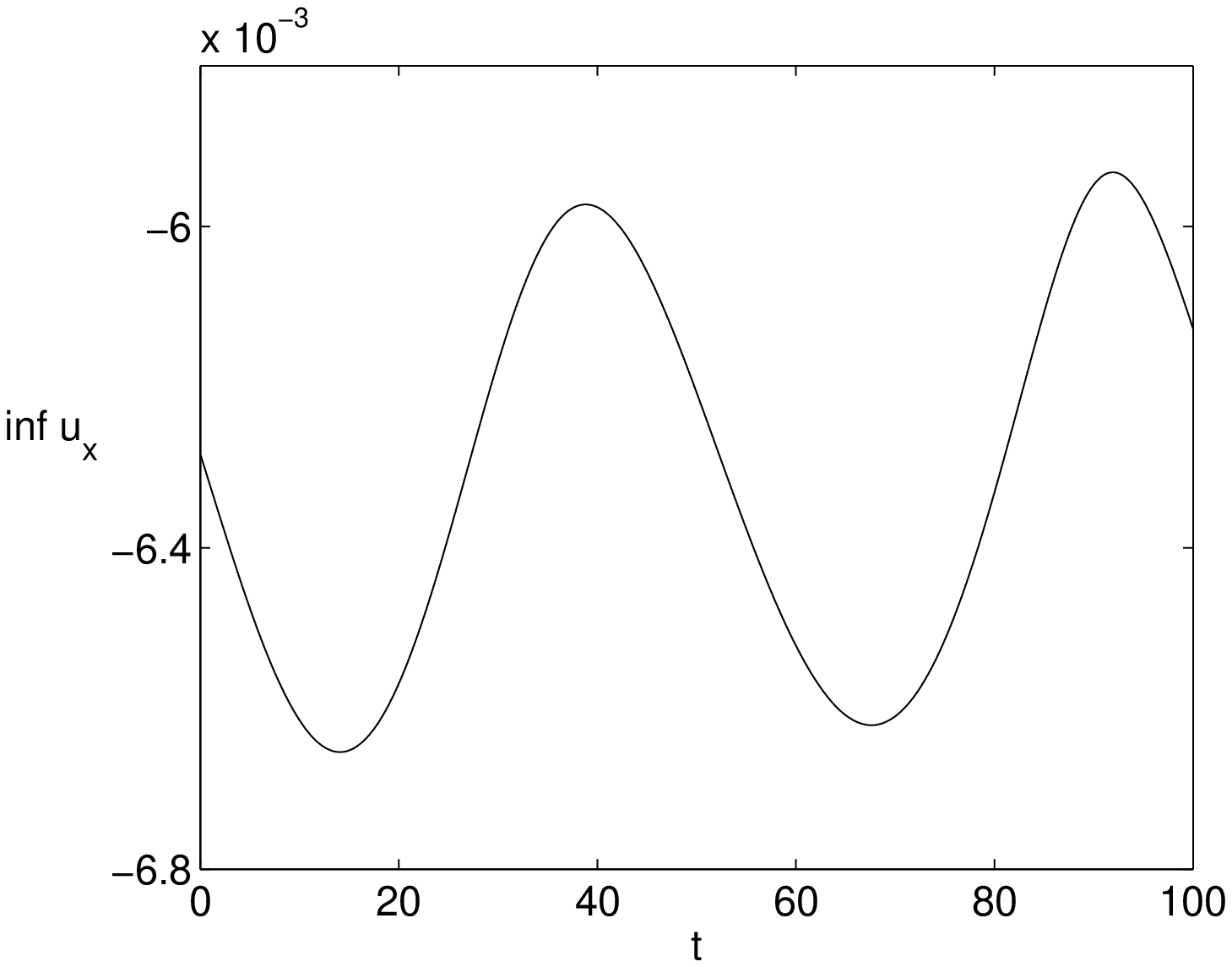}\\
\includegraphics[width=0.45\textwidth]{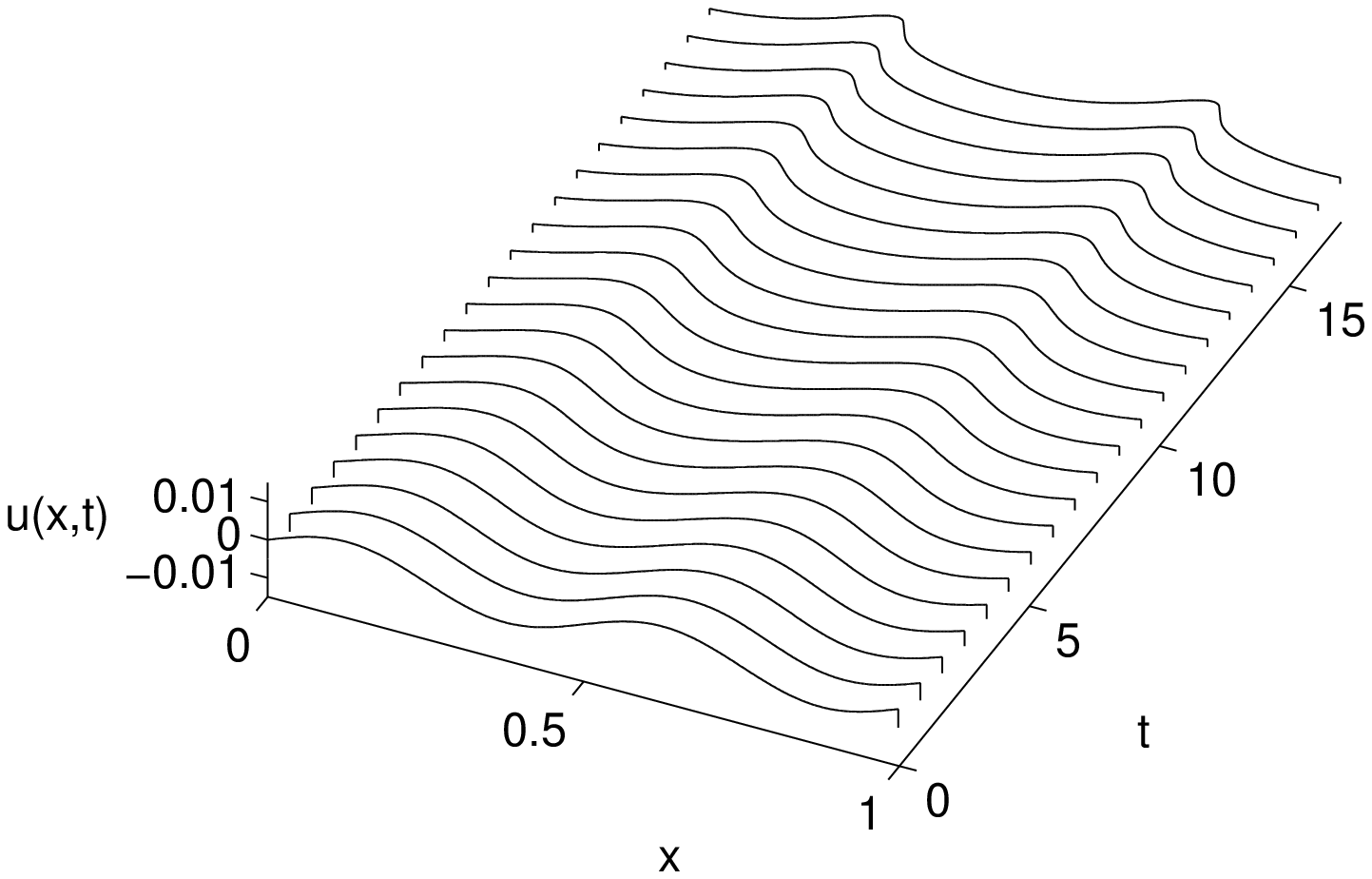} & \includegraphics[width=0.45\textwidth]{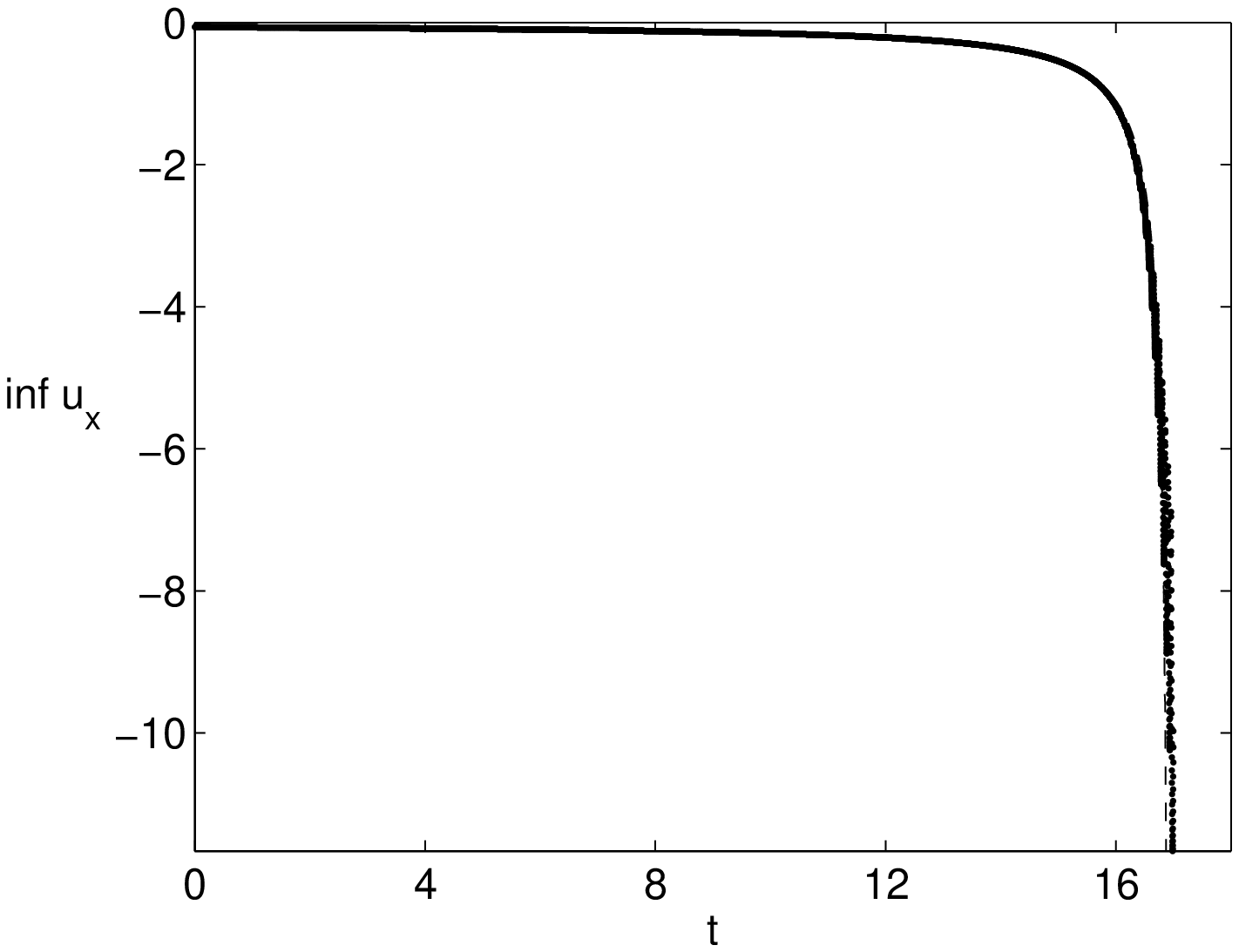}\\
\end{tabular}
\end{center}
\caption{The same as Figure \ref{fig-1} but for $a=0$, $b=0.0005$
(top) and $a=0$, $b=0.005$ (bottom). The least squares fit  is
computed with $C \approx -1.060$ and $B \approx 16.964$.}
\label{fig-3}
\end{figure}

Numerical simulations of the Ostrovsky-Hunter equation
(\ref{OstHunter}) for initial data (\ref{initial-data}) are
performed with the pseudo--spectral method for $N = 4096$ Fourier
harmonics with the time step of $dt = 0.001$. Figures
\ref{fig-1},\ref{fig-2}, and \ref{fig-3} show two dynamical
evolutions for three cases  $b = 0$, $a = 2b$, and $a = 0$. In all
cases, no wave breaking occur for sufficiently small values of
$(a,b)$ (far below the lower bound on Figure
\ref{fig-parameter-plane}) but the wave breaking does occur if the
values of $(a,b)$ are selected to be larger (still below the lower
bound on Figure \ref{fig-parameter-plane}). Thus, we conclude that
none of the wave breaking criteria is sharp.

Right panels of Figures \ref{fig-1},\ref{fig-2}, and \ref{fig-3}
show the behavior of $\inf_{x \in \mathbb{S}} u_x(t,x)$ versus
$t$. When the wave breaking occurs (bottom panels of each figure),
we compute the linear regression $B + Ct$ of
$-(\inf_{x\in\mathbb{S}} u_x(t,x))^{-1}$, where $(B,C)$ are
constants. According to Theorem \ref{theorem-blowup-rate}, $B
\approx T$ and $C \approx -1$ near the singularity, so that $B$
can be taken as an approximation for the blow-up time $T$ and $|C
+ 1|$ can be taken as an estimate for the error of the linear
regression. The numerical values on Figures
\ref{fig-1},\ref{fig-2}, and \ref{fig-3} show that $C$ is close to
$-1$ by the errors in $1\%$, $4\%$, and $6\%$, respectively.

\begin{figure}
\begin{center}
\begin{tabular}{cc}
\includegraphics[width=0.45\textwidth]{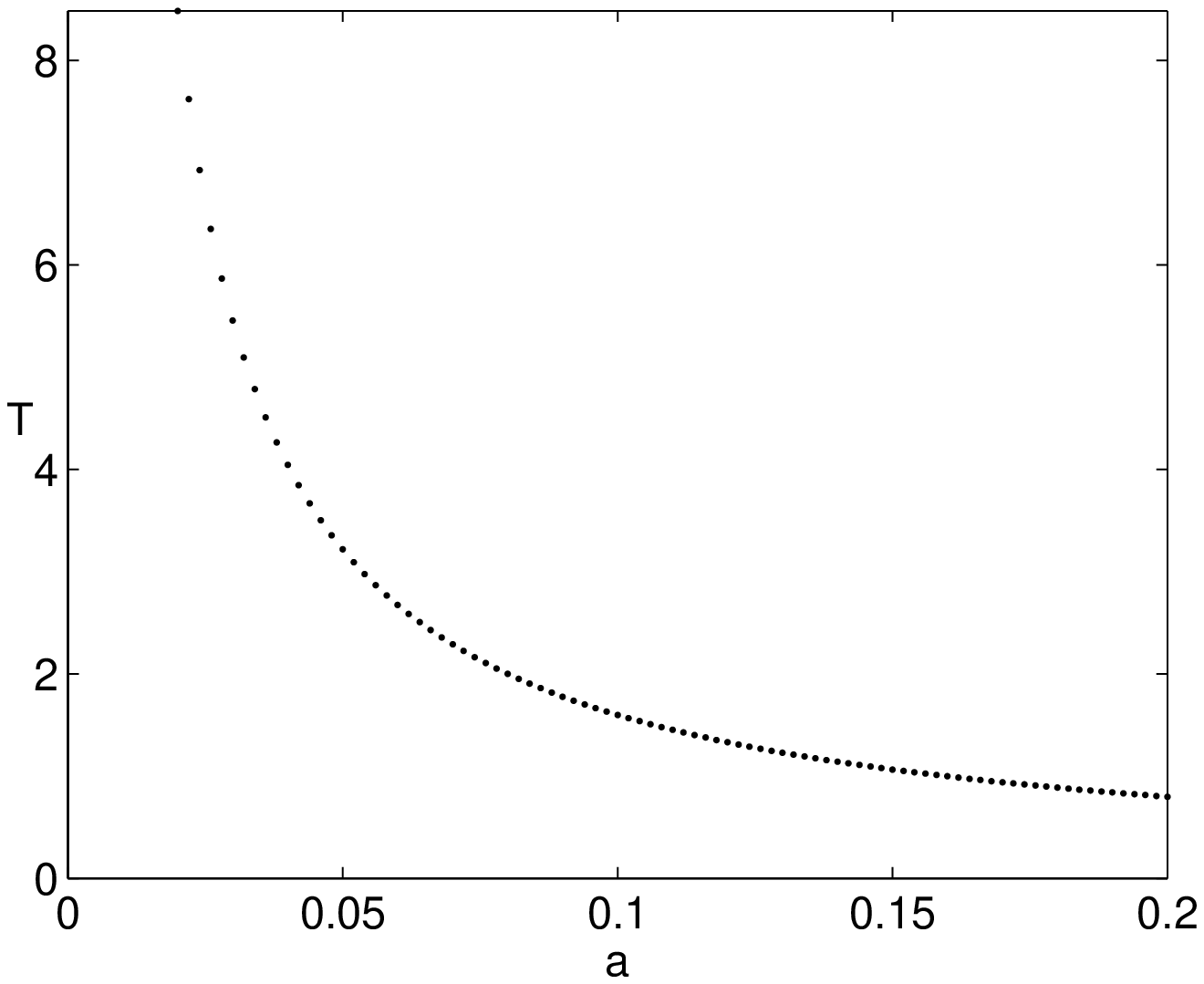} &
\includegraphics[width=0.45\textwidth]{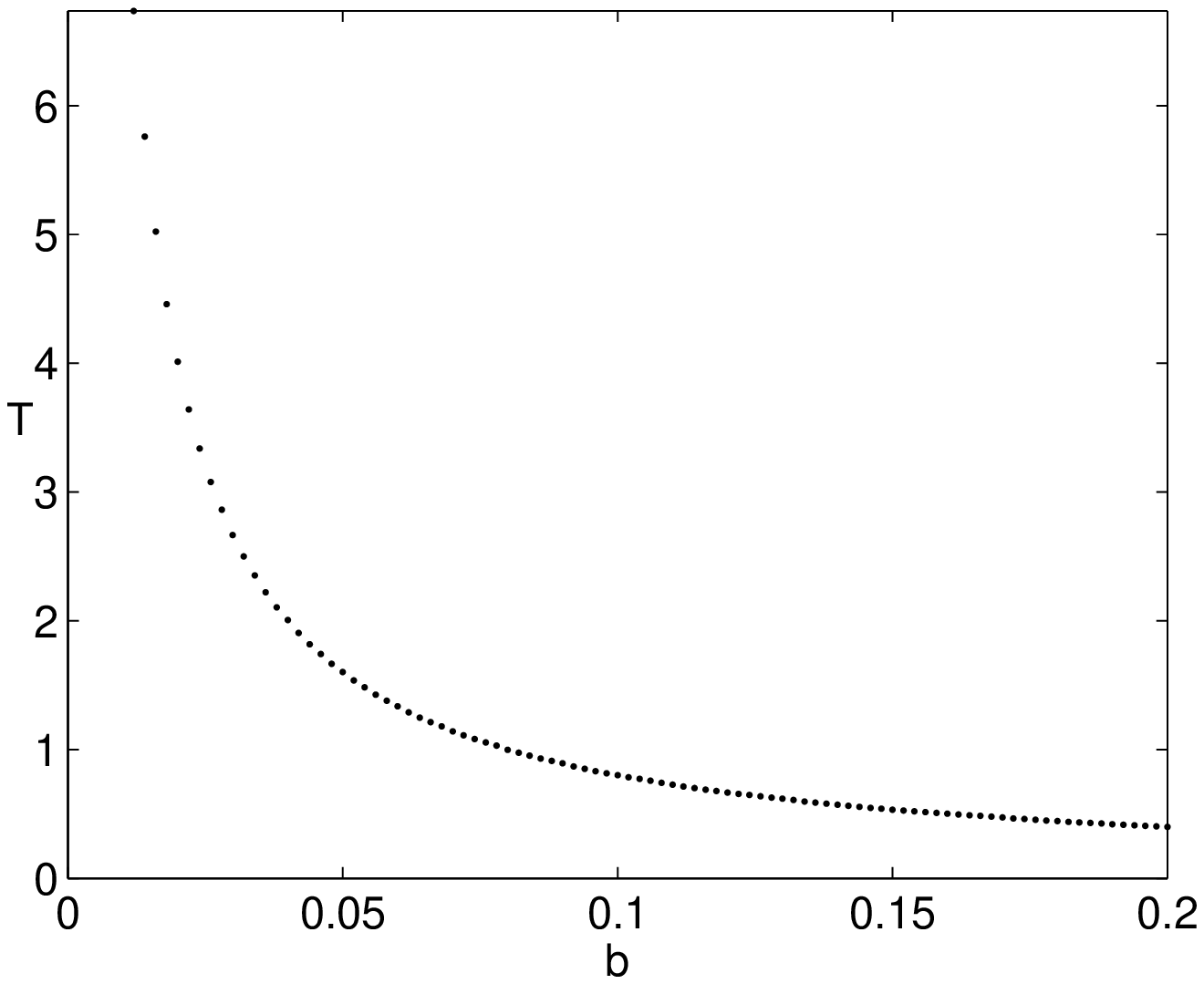}
\end{tabular}
\end{center}
\caption{Estimates of the blow-up time $T$ versus $a$ for $b = 0$
(left) and $b$ for $a = 0$ (right).} \label{fig-blowup}
\end{figure}

Finally, Figure \ref{fig-blowup} shows the blow-up time $T$
estimated by the above technique versus parameters $a$ for $b = 0$
and parameter $b$ for $a = 0$. We can clearly see that the wave
breaking holds for $(a,b)$ in the shaded area of Figure
\ref{fig-parameter-plane}. The blow-up time $T$ becomes smaller
for larger values of $(a,b)$.

\appendix
\section{Appendix A: \\Nonexistence of localized traveling-wave solutions}

Consider the differential equation
\begin{equation}
\label{traveling}
(c-\varphi(x)) \varphi''(x) =(\varphi'(x))^2 -  \gamma \varphi(x),  \quad x \in \mathbb{R}
\end{equation}
for traveling-wave solutions $u(t,x) = \varphi(x-ct)$ of the Ostrovsky--Hunter equation
(\ref{OstHunter}), where $c \in \mathbb{R}$ is the wave speed.

\begin{theorem}
There are no nontrivial solutions of (\ref{traveling}) with $c \in \mathbb{R}$ such that
$\varphi \in C^2(\mathbb{R})$ and $\lim_{|x| \to \infty} \varphi(x) = \lim_{|x| \to \infty} \varphi'(x) =0$.
\end{theorem}

\begin{proof}
Arguing by contradiction, we assume the existence of $\varphi \in C^2(\mathbb{R})$ with
$\lim_{|x| \to \infty} \varphi(x) = \lim_{|x| \to \infty} \varphi'(x) = 0$ as a solution of (\ref{traveling}). Let $v(x)$ be defined by
\begin{equation}
\label{3-3}
 - c \varphi + \frac{1}{2} \varphi^2 = \gamma v,
\end{equation}
so that $\varphi(x) = v''(x)$. Multiplying (\ref{3-3}) by
$\varphi'(x)$ and taking integral over the interval $ ( - \infty, \, x) $ yields
\begin{equation}
\label{3-4}
-\frac {c} {2} \varphi^2  + \frac {1}{6} \varphi^3 + \frac {\gamma}{2} (v')^2 - \gamma \varphi v = 0.
\end{equation}
By equation (\ref{traveling}), we have
$\int_{\mathbb{R}} \varphi (x) \, dx = 0$. Since $\varphi \in C^2(\mathbb{R})$, there exists a smallest zero point $ \xi_1 \in \mathbb{R}$ for $ \varphi(x)$   in the sense of $ \varphi (\xi_1) = 0 $ and $ \varphi(x) \ne 0$ for $ x \in ( - \infty, \, \xi_1)$.  On the other hand, it is deduced from (\ref{3-4}) that $v'(\xi_1) = 0$ and
$\lim_{x \to - \infty} v'(x) = 0$. By Rolle's theorem, there exists a point $\xi_2 \in ( - \infty, \, \xi_1) $ such that $ \varphi(\xi_2) = v''(\xi_2) = 0$, which contradicts the assumption on the smallest $ \xi_1 $ with $ \varphi(\xi_1) = 0. $  This completes the proof of the theorem.
\end{proof}

\end{document}